\documentclass[]{article}
\usepackage{graphicx}
\usepackage{graphics,color}
\usepackage{multicol}
\usepackage{multirow}
\usepackage{amssymb,amsmath,amsthm,epsfig,bbm}
\usepackage{dsfont}

\newcommand{\N}{\mathbb{N}}
\newcommand{\Z}{\mathbb{Z}}

\newcommand{\pu}{\pmb{u}}
\newcommand{\pw}{\pmb{w}}
\newcommand{\pk}{\pmb{k}}
\newcommand{\px}{\pmb{x}}
\newcommand{\B}{\mbox{B}}

\newtheorem{theorem}{Theorem}
\newtheorem{lemma}{Lemma}
\newtheorem{corollary}{Corollary}
\newtheorem{definition}{Definition}
\newtheorem{proposition}{Proposition}

\newtheorem{remark}{Remark}
\newcommand{\R}{\mathbb R}

\usepackage{color}

\newcounter{reh}
\newcounter{rek}
\newcommand{\PP}{\mathbb{P}}

\begin{document}

\begin{center}
{\large {\bf {A Stable Jacobi polynomials based least squares  regression estimator associated with 
an ANOVA decomposition model.}}}\\
\vskip 0.5cm Mohamed Jebalia$^a$ and Abderrazek Karoui$^b$ 
\footnote{
\noindent
Emails: abderrazek.karoui@fsb.ucar.tn  (A. Karoui, corresponding author), mohamed.jebalia@enib.ucar.tn (M. Jebalia).\\
This work was supported in part by the  
DGRST  research grant LR21ES10 and the PHC-Utique research project 20G1503.}  
\end{center}

\vskip 0.5cm
\noindent
$^a$  University of Carthage, National School of Engineering of Bizerte, Menzel Abderrahman 7035,  Tunisia.\\
$^b$ University of Carthage, Faculty of Sciences of Bizerte, Jarzouna 7021,  Tunisia.\\


\noindent{\bf Abstract}--- 
In this work, we construct a stable and fairly fast estimator for solving non-parametric multidimensional regression problems. The proposed estimator is based on the use of multivariate Jacobi polynomials that generate a basis for  a reduced size
of $d-$variate finite dimensional  polynomial space. An ANOVA decomposition trick has been used for building this later polynomial space. Also, by using some results from the theory of positive definite random matrices, we show that the proposed estimator is stable under the condition that the i.i.d. random sampling points for the different  covariates of the regression problem, follow a $d-$dimensional Beta distribution. Also, we provide the reader with an estimate for the $L^2-$risk error of the estimator. Moreover, a more precise estimate of the 
quality of the approximation
is provided under the condition that the regression function belongs to some weighted Sobolev space. Finally, 
the various theoretical results of this work
are supported by numerical simulations. \\


\section{Introduction} In this work, we combine the popular technique of generalized polynomials chaos (gPC)
\cite{Tarakanov, Xiu} and a special family of $d-$variate Jacobi polynomials 
in order
to solve a $d-$dimensional non-parametric regression problem. This last problem is one of the important as well as an active research topic from the machine learning area, see for example \cite{Tarakanov, Torre}. Note that a machine learning algorithm can be briefly described as an algorithm for the approximation of an unknown function $f$ that maps in general a random vectors $X \in \mathbb R^d$ to an observed real valued variable $Y.$ An estimator or an approximation $\widehat  f$ of $f$ is constructed by the use a training data set $\{ (X_i,Y_i),\, 1\leq i\leq n\}.$ Note that unlike a parametric learning algorithm, 
where $\widehat f$ is given in terms of a set of fixed size of parameters, a non-parametric learning algorithm does not require any assumption about the function $f$ or its estimator $\widehat f.$ Usually, for a non-parametric (NP) model, the function $f$ lies in an infinite dimensional functional space. Consequently, the NP models have the advantage to better fit a wide range of the true functions $f$. We should mention that the   multidimensional NP regression problem is  frequently encountered in a wide range of scientific fields. An NP  learning algorithm for solving this problem aims to provide a convenient estimate $\widehat f$ for the true regression function $f$ associated  with the regression problem,
\begin{equation}
	\label{eq:Pb1}
	Y_i=f(X_i)+\varepsilon_i, \  i=1,\ldots,n \,.
\end{equation}
Here,  the $X_i \in \mathbb R^d$ are assumed to be  i.i.d. random  vectors and the 
 $(\varepsilon_i)_{1\le i \le n}$ are centered i.i.d. random variables with a finite variance 
$\sigma^2=\mathbb E[\varepsilon_1^2]$ and independent from the $X_i.$  The goal of a learning algorithm is to minimize the empirical risk $\mathcal R_{emp}(f)$ over a given class of functional space $\mathcal H.$ That is to solve the minimization problem 
\begin{equation}\label{minimization_pb}
\widehat f =\arg\min_{f\in \mathcal H} \mathcal R_{emp}(f) =\arg\min_{f\in \mathcal H} \frac{1}{n} \sum_{i=1}^n  L(f(X_i), Y_i),
   \end{equation}
 where, $L(\cdot,\cdot)$ is a non-negative loss function. Among the frequently used loss functions from the literature, we cite the Tikhonov regularized loss \cite{smale}  and the weighted $\ell_2$ loss function \cite{Shin1}, given respectively by 
 $$L_\lambda(F(X_i), Y_i)=\big(f(X_i)-Y_i\big)^2+\lambda \|f\|^2_{\mathcal H},\qquad L_\omega (f(X_i),Y_i)= \omega_i \big(f(X_i)-Y_i\big)^2,$$
 for some convenient regularization parameter $\lambda >0$ and finite weight sequence $(\omega_i)_i.$ In general, the least squares scheme is used to solve the minimization problem \eqref{minimization_pb}. \\
 
 It is well known that despite the superiority of the NP model in terms of quality of approximation of the true functional $f$, it suffers from the curse of dimensionality for large or even moderate values of $d,$ the number of covariates. The computational load by an NP algorithm grows fast with the dimension $d$ and its convergence rate or its associated risk error rate slows drastically. For instance, it has been shown in \cite{Stone}, see also \cite{Bauer, Gyorfi} that if $f$ is of class  $C^p$ with $p-$th derivative being H\"older continuous, then the optimal convergence rate of any NP least-squares  estimator $\widehat f_n,$ solution of the minimization problem \eqref{minimization_pb} is given by 
  ${\displaystyle \mathbb E\Big[\| f-\widehat f_n \|_{\mathcal H}^2\Big]=O\Big(n^{-2p/(2 p+d)}\Big).}$ In order to reduce the computational load required by an algorithm for solving more general models with random inputs or models from uncertainty quantification (UQ) area, a popular technique of polynomial chaos expansion (PCE) is successfully used in the literature. This technique aims to approximate the output variable $Y$ by using a projection over a reduced size of orthogonal polynomial basis. The PCE scheme  has been first introduced by N. Wiener in his pioneer work \cite{Wiener}, for the Hermite polynomials and Gaussian random variables.
  Recently, there is a growing interest in the study and the use for UQ applications of a generalized version of PCE, called generalized polynomial chaos (gPC). The gPC was first introduced by \cite{Xiu} and it aims to extend the PCE to various discrete and continuous probability distributions associated with the set of weight functions for the family of orthogonal polynomials of the Askey--scheme. For more details on PCE and gPC schemes and their associated UQ related applications, the reader is refereed to \cite{Hadigol, Hampton1, Jakeman1,  Jakeman2, Loukrezis, Torre, Zhou2}.  Nonetheless, the gPC based learning algorithm  still has the  limitation to be slow for moderate large values of the dimension $d$. To overcome this problem,  various solutions have been considered in the literature. Among the popular adopted solutions, we cite the use of sparsity and optimal sampling techniques, see for example \cite{Blatman, Castro, Guo, Hampton2, Jakeman2, Lin, Pan, Shin2}, dimension reduction through a sensitivity analysis techniques, \cite{Alexanderian, Blatman2},  as well as the use   of  partial functional ANOVA decomposition technique, see for example \cite{Huang, Lin2, Potts, Saltelli, Sobol, Stone2}. 
  More precisely, a partial  functional ANOVA decomposition consists 
  in the approximation of a real valued function $f(x_1,\cdots,x_d)$ by a sum of functions with reduced number of variables of the form $f_i(x_{i_1},\ldots,x_{i_m})$ where $1\leq m \leq  d$ and the $i_j\in [[1,d]].$ Note that 
  for
  the full ANOVA decomposition, that is $m=d,$ the uniqueness of the ANOVA decomposition of a function $f\in L^1(J^d),\, J=[0,1]$ has been shown in \cite{Sobol}. Also, among the popular reduced size multivariate polynomials spaces used by a gPC scheme, we cite  the total degree  space of degree $N,$ given by 
  $\mathcal P_N^{TD}=\mbox{Span}\{\pmb x^{\pmb i}=x_1^{i_1}\cdots x_d^{i_d},\quad \|\pmb i\|_1=\sum_{j=1}^d |i_j| \leq N \}$
  and the hyperbolic cross space of  degree $N,$ given by 
  $\mathcal P_{q,N}^{HC}=\mbox{Span}\{\pmb x^{\pmb i}=x_1^{i_1}\cdots x_d^{i_d},\quad \|\pmb i\|_q=\Big(\sum_{j=1}^d |i_j|^q \Big)^{1/q}\leq N \},$ $0<q<1.$ For more details on these polynomials spaces, the reader is refereed to \cite{Blatman}. \\
  
  In this work, we introduce a new $d-$variate polynomial space constructed from uni-variate Jacobi polynomials associated with parameters $\beta=\alpha \geq -\frac{1}{2}$ and orthonormal over$I^d=[-1,1]^d.$ More precisely, for two  positive integers $N\geq 1$ and 
  $1\leq m \leq \min(d,N),$ we let the ANOVA type Jacobi polynomials space 
  \begin{equation}
	\label{eq:Pspace1}
\mathcal{P}_{N,m,d}= Span\left\{ \prod_{i \in \pmb u} \widetilde{P}_{k_i}^{(\alpha)}(x_i) \ ; \ \pu \subset [[1,d]] \ ;
\ \ |\pu| \le m ; \ \pk \in \mathbb N_0^d, \ \ ||\pk||_1 \le N  \right\}\,.
\end{equation}
Here, the $\widetilde{P}_{k_i}^{(\alpha)}$ are the orthonormal  uni-variate Jacobi polynomials , associated with a parameter $\alpha\geq -\frac{1}{2}.$ The dimension of  
 $\mathcal{P}_{N,m,d}$ is given by   ${\displaystyle  M_{N,m,d}=\dim \mathcal{P}_{N,m,d}= \sum_{k=0}^m {d \choose k}{N \choose k}.}$ For the extreme case $m=\min(d,N),$ the space $\mathcal{P}_{N,m,d}$ is reduced to the   total degree  space polynomials of degree $N,$ given by  $\mathcal P_N^{TD}.$ One of the main results of this work is to prove the stability of  our proposed  least-squares estimator 
 $\widehat{f}_{N,n,m}^{(\alpha)}$ which is the solution of the minimization problem \eqref{minimization_pb} with the functional space $\mathcal H=\mathcal{P}_{N,m,d} $ and a loss function $L$  given by the Euclidean distance of $\mathbb R^n.$ Note that there is a growing interest in the study of the stability issue of estimators of functions with random inputs, see for example \cite{Adcock, Cohen1, Cohen, Migliorati, Narayan}. In the present work,  we show that under the condition that the $X_i$ follow a $d-$variate $\mbox{Beta}(\alpha+1,\alpha+1)$ distribution with support $I^d,$ the positive definite  $n\times n$ random matrix involved in the construction of our proposed least-squares based estimator $\widehat{f}_{N,n,m}^{(\alpha)}$ is with high probability well conditioned in the $2-$norm. Moreover, we give an estimate for the $L^2-$risk error of a truncated version of the $\widehat{f}_{N,n,m}^{(\alpha)}$ which we denote by $\widehat{F}_{N,n,m}^{(\alpha)}.$ More precisely, if $\|\cdot\|_\alpha$ denotes the weighted $L^2-$norm associated with the weight $\pmb \omega_\alpha=\prod_{i=1}^d (1-x_i^2)^\alpha,$
then we give an estimate of the $L^2-$risk error  $\mathbb E \Big[ \|f-
\widehat{F}_{N,n,m}^{(\alpha)}\|^2_\alpha\Big].$ This later is given in terms of the classical bias-variance decomposition. The variance term decays at a rate of $O\big(\frac{M_{N,m,d}}{n}\big).$ Here, $M_{N,m,d}$ is the dimension of our proposed polynomial space $\mathcal P_{N,m,d}.$ The bias term f the $L^2-$risk involves the quantity $\|f- \Pi_{N,m} f\|_\alpha,$ where $\Pi_{N,m}$ denotes the orthogonal projection over  $\mathcal P_{N,m,d}.$ An estimate of this last quantity is given under the hypothesis that the true regression function lies in a weighted Sobolev space with given Sobolev smoothness property.\\

This work is organized as follows. In section 2, we give some mathematical preliminaries on Matrix Chernoff eigenvalues bounds, the Gershgorin circle theorem for bounding the spectrum of a square matrix, as well as some properties of the Jacobi polynomials. In particular, we give some  bounds of these polynomials that will be used for proving different results of this work. In section 3, we describe the new adopted and reduced size  ANOVA type multivariate Jacobi polynomial space $\mathcal P_{N,m,d}.$ Moreover, we  provide the reader with an estimate of the  dimension of this later. Section 4 of this work is devoted to the proof of the stability of the  proposed least-squares and gPC based NP regression estimator $\widehat{f}_{N,n,m}^{(\alpha)}.$ In section 5, we give an estimate for the weighted $L^2-$risk error of a truncation version  of the estimator $\widehat{f}_{N,n,m}^{(\alpha)}.$ Moreover, in section 6, we give an estimate for the bias term of the previous weighted $L^2-$risk error, when the true regression function belongs to some weighted Sobolev space. Finally, in section 7, we give some numerical simulations that illustrate the different results of this work.\\

\section{Mathematical Preliminaries and estimates for Jacobi polynomials}
\label{Sec:Prelim}

In this paragraph, we first give some mathematical preliminaries from the literature that will be used frequently in this work. Then, we give some useful estimates for the Jacobi polynomials. These estimates are needed for the proof of the stability property of our proposed special Jacobi polynomials  multivariate non-parametric (NP) regression estimator.

\subsection{Mathematical preliminaries}

We first recall the {\bf Matrix Chernoff Theorem} (see for example \cite{Tropp}) and the {\bf Gershgorin circle Theorem} (see for example \cite{Horn}) that will be useful 
 for proving the stability of our NP regression estimator.\\
 
{\bf Matrix Chernoff Theorem:}
Consider a sequence of $n$ independent $D \times D$
random Hermitian matrices $\{\pmb{Z}_k\}$. Assume that for some $L>0$, we have 
$$ 0 \preccurlyeq \pmb{Z}_k \preccurlyeq L. \pmb{I}_D\,.$$
Let 
$$ \pmb{A}=\sum_{k=1}^n \pmb{Z}_k, \ \mu_{\min}=
\lambda_{\min}(\mathbb{E}(\pmb{A})). \  
 \mu_{\max}=
\lambda_{\max}(\mathbb{E}(\pmb{A}))\,.
$$
Then, for any $\delta \in (0,1]$, we have 
\begin{equation}
    \PP\left(\lambda_{\min}(\pmb{A} \le (1-\delta) \mu_{\min} \right) \le D. \exp\left(  -\frac{\delta^2 \mu_{\min}}{2L}\right), \ 
    \PP\left(\lambda_{\max}(\pmb{A} \ge (1+\delta) \mu_{\max} \right) \le D. \exp\left(  -\frac{\delta^2 \mu_{\max}}{3L}\right)
\end{equation}
\\
{\bf Gershgorin circle Theorem:} Let $A=[a_{i,j}]_{1\le i,j \le n}$ be a complex matrix.
For $1 \le i \le n$, let $R_i=\sum_{j \neq i} |a_{i,j}|$.
Then every eigenvalue of  $A$ lies within at least one of the discs $D(a_{ii},R_i)$.\\~\\
In the sequel, we let $\Gamma(a)$ and $\B(a,b)$ respectively denote the usual Gamma and  Beta functions with $a,b>0$.\\
For an integer $k\ge 0$ and  $\alpha \ge -\frac{1}{2}$, let $\widetilde{P}_k^{(\alpha,\alpha)}$ denote the normalized Jacobi polynomial defined on $I=[-1,1]$ of degree $k$ and parameters $(\alpha,\alpha)$. We have: 
\begin{equation}\label{JacobiP}
\widetilde{P}_k^{(\alpha,\alpha)}(x)=\frac{1}{\sqrt{h_k^{(\alpha)}}}P_k^{(\alpha)}(x) ,
\  h_k^{(\alpha)}=
\frac{2^{2\alpha+1}\Gamma^2(k+\alpha+1)}{k!(2k+2\alpha+1)\Gamma(k+2\alpha+1)}
\,.
\end{equation}
The polynomials $\widetilde{P}_k^{(\alpha,\alpha)}, \ k \ge 0$
satisfy the orthonormality relation
$$\int_I \widetilde{P}_j^{(\alpha,\alpha)} \widetilde{P}_k^{(\alpha,\alpha)} w_(\alpha)(x)dx=\delta_{j,k},
\, \ w_{\alpha}(x)=(1-x^2)^{\alpha}\,.$$
In the sequel, in order to alleviate notations, we will use the notation $\widetilde{P}_k^{(\alpha)}$
instead of $\widetilde{P}_k^{(\alpha,\alpha)}.$

The following Lemma regroups different useful identities and inequalities that can be easily found in the literature, see for example 
\cite{Andrews, Batir, Karoui-Souabni}. 
\begin{lemma} \label{lm:useful}
Let  $J_{a}$ be the Bessel function of the first kind and order $a>-1$. Then, we have 
\begin{enumerate}
    \item For any $x \in \R$ and for any integer $m\geq 0,$ 
    \begin{equation}
    \label{eq:Bessel1}
        \int_{-1}^{1} e^{ixy}\widetilde{P}_m^{(\alpha)}(y)w_{\alpha}(y)dy=
        i^m\sqrt{\pi}.\sqrt{2m+2\alpha+1}
        \sqrt{\frac{\Gamma(m+2\alpha+1)}{\Gamma(m+1)}}
        \frac{J_{m+\alpha+1/2}(x)}{x^{\alpha+1/2}}\,.
    \end{equation}
    \item For any  $x \in \R$ and any real  $\mu >-1$, we have 
    \begin{equation}
    \label{eq:Bessel2}
        |J_{\mu}(x)| \le \frac{|x|^{\mu}}{2^{\mu}\Gamma(\mu+1)}\,.
    \end{equation}
\item For $ x >-\frac{1}{2}$, 
\begin{equation}
	\label{eq:Gamma}
\sqrt{2e} \left(\frac{x+\frac{1}{2}}{e}\right)^{x+\frac{1}{2}} \le \Gamma(x+1)
\le \sqrt{2\pi} 
\left(\frac{x+\frac{1}{2}}{e}\right)^{x+\frac{1}{2}}\,.
\end{equation}
\end{enumerate}
\end{lemma}

\subsection{Estimates for Jacobi polynomials}
In order to provide estimates for Jacobi polynomials, we will need the following result.
\begin{lemma}
\label{lm:h}
   For $\alpha \geq -\frac{1}{2},$ the function $h_0^{(\alpha)}=2^{2\alpha+1}
    \text{Beta}(\alpha+1,\alpha+1) $ is bounded 
    as follows
    \begin{equation}
    \frac{1}{C^2(\alpha)} \le    h_0^{(\alpha)} \le 2\pi\,,
    \end{equation}
    where  $C(\alpha):=\left[\frac{\pi}{e^2}
        \left(\alpha+\frac{3}{4}\right)\right]^{\frac{1}{4}}$.
\end{lemma}
\proof:
We first establish the lower bound of $h_0^{(\alpha)}$ using (\ref{eq:Gamma}).
\begin{eqnarray*}
\frac{1}{h_0^{(\alpha)}}
&=&\frac{1}{2^{2\alpha+1}}
	\frac{\Gamma(2\alpha+2)}{\Gamma^2(\alpha+1)}\\
	&\leq & 
\sqrt{\pi} {e} ^{-\frac{3}{2}}
	\left(\alpha+\frac{3}{4} \right)^{2\alpha+\frac{3}{2}}
	\left(\frac{1}{\alpha+\frac{1}{2}} \right)^{2\alpha+1}\\
&\leq &
	\sqrt{\pi} {e} ^{-\frac{3}{2}}
		\left(\alpha+\frac{3}{4} \right)^{\frac{1}{2}}
	\left(1+\frac{1}{2}\frac{1}{2\alpha+1} \right)^{2\alpha+1}
	\le 
	\sqrt{\pi} {e} ^{-1}
		\left(\alpha+\frac{3}{4} \right)^{\frac{1}{2}}.	
\end{eqnarray*}
For the upper bound of $h_0^{(\alpha)}$, we will use again (\ref{eq:Gamma}). We get
\begin{eqnarray*}
h_0^{(\alpha)}
=2^{2\alpha+1}
	\frac{\Gamma^2(\alpha+1)}{\Gamma(2\alpha+2)}
&	\leq&   2^{2\alpha+2} \pi 
		\left(\frac{\alpha+\frac{1}{2}}{e} \right)^{2\alpha+1}
		\frac{1}{\sqrt{2e}}
			\left(\frac{e}{2\alpha+\frac{3}{2}} \right)^{2\alpha+\frac{3}{2}}
\\
& = &
2^{2\alpha+\frac{3}{2}} \pi 
		\left(\frac{\alpha+\frac{1}{2}}{e} \right)^{2\alpha+1}
		\frac{1}{\sqrt{e}}
			\left(\frac{e}{2\alpha+\frac{3}{2}} \right)^{2\alpha+\frac{3}{2}}
			\\
			& = &
 \pi 
			\left(\frac{1}{\alpha+\frac{3}{4}} \right)^{\frac{1}{2}}
		\left(\frac{\alpha+\frac{1}{2}}{\alpha+\frac{3}{4}} \right)^{2\alpha+1}
						 \leq 2\pi.
\end{eqnarray*}
\qed~\\

The following proposition provides us with some useful estimates for the normalized Jacobi polynomials $\widetilde{P}_k^{(\alpha)}.$

\begin{proposition}[Bounds for $||\widetilde{P}_k^{(\alpha)}||_{\infty}$]\label{prop:boundsalpha}~ Under the previous notation, let $\alpha\geq -\frac{1}{2},$ then for any integer $k\geq 1,$ we have
\begin{equation}
\label{Jacobi_ineq1}
\big\|\widetilde{P}_k^{(\alpha)}\big\|_{\infty}
        \leq \left\{\begin{array}{ll} 
       \frac{\pi^{\frac{1}{4}} e^{\alpha+\frac{3}{4}}}{\sqrt{2}} k^{\alpha+\frac{1}{2}}
& \mbox{ if } \alpha > -\frac{1}{2}\\ 
~ & ~\\ 
\frac{2}{\sqrt{\pi}} &\mbox{ if } \alpha=-\frac{1}{2}.
\end{array}\right.
\end{equation}
Moreover, for $k=0,$ 
\begin{equation}\label{Jacobi_ineq2}
\big\|\widetilde{P}_0^{(\alpha)}\big\|_{\infty}
        \le \ C(\alpha),
        \quad \alpha\geq -\frac{1}{2}\,,
\end{equation}
where $C(\alpha)$ is the quantity defined in 
Lemma~\ref{lm:h}.
    \end{proposition}
    
\proof: 
Let $k \ge 1$, since  
$$ \frac{1}{\sqrt{h_k^{(\alpha)}}}=\frac{1}{2^{\alpha}
\Gamma(k+\alpha+1)}
\sqrt{\Gamma(k+1)\Gamma(k+2\alpha+1)} \sqrt{k+\alpha+\frac{1}{2}},$$
then 
$$ \frac{||P_k^{(\alpha)}||_{\infty}}{\sqrt{h_k^{(\alpha)}}\sqrt{k+\alpha+\frac{1}{2}}}= \frac{1}{2^{\alpha}\Gamma(\alpha+1)}\sqrt{\frac{\Gamma(k+2\alpha+1)}{\Gamma(k+1)}}\,.  $$
For  $\alpha=-\frac{1}{2}$, one gets
$$ \|\widetilde{P}_k^{(-\frac{1}{2})}||_{\infty}=
\frac{||P_k^{(-\frac{1}{2})}\|_{\infty}}{\sqrt{h_k^{(1/2)}}}=\frac{\sqrt{2}}{\Gamma(1/2)}\sqrt{\frac{\Gamma(k)}{\Gamma(k+1)}}
\sqrt{k+1}=\sqrt{\frac{2}{\pi }}\sqrt{1+\frac{1}{k}}
\le \frac{2}{\sqrt{\pi} }.$$ 
Next, for the case  $\alpha >-\frac{1}{2}$ and by 
Using the bounds of the Gamma function~(\ref{eq:Gamma}), one gets 
$$ \frac{\Gamma(k+2\alpha+1)}{\Gamma(k+1)} 
\le \left( \frac{\pi}{e}\right)^{\frac{1}{2}}
\left( \frac{k}{e}\right)^{2\alpha}
\left(1+\frac{2\alpha}{k+\frac{1}{2}} \right)^{k+\frac{1}{2}}
\left(1+\frac{2\alpha+\frac{1}{2}}{k} \right)^{2\alpha}\,.
$$
Thus, 
$$ \frac{||P_k^{(\alpha)}||_{\infty}}{\sqrt{h_k^{(\alpha)}}\sqrt{k+\alpha+\frac{1}{2}}} \le  \frac{1}{2^{\alpha}\Gamma(\alpha+1)}
\left( \frac{\pi}{e}\right)^{\frac{1}{4}}
\left( \frac{k}{e}\right)^{\alpha}
\left(1+\frac{2\alpha}{k+\frac{1}{2}} \right)^{\frac{k+\frac{1}{2}}{2}}
\left(1+\frac{2\alpha+\frac{1}{2}}{k} \right)^{\alpha}
\,.  $$
Since $\left(1+\frac{2\alpha}{k+\frac{1}{2}} \right)^{\frac{k+\frac{1}{2}}{2}} \le e^{\alpha} $, then we get
$$ \frac{||P_k^{(\alpha)}||_{\infty}}{\sqrt{h_k^{(\alpha)}}} \le  \frac{k^{\alpha+\frac{1}{2}}}{2^{\alpha}\Gamma(\alpha+1)}
\left( \frac{\pi}{e}\right)^{\frac{1}{4}}
\left(2\alpha+\frac{3}{2} \right)^{\alpha} 
\sqrt{\alpha+\frac{3}{2}}
=
\frac{k^{\alpha+\frac{1}{2}}}{\Gamma(\alpha+1)}
\left( \frac{\pi}{e}\right)^{\frac{1}{4}}
\left(\alpha+\frac{3}{4} \right)^{\alpha} 
\sqrt{\alpha+\frac{3}{2}}
\,.  $$
Using again~(\ref{eq:Gamma}), we get 
$\frac{1}{\Gamma(\alpha+1)} \le \frac{1}{\sqrt{2e}}
\left(\frac{e}{\alpha+\frac{1}{2}} \right)^{\alpha+
\frac{1}{2}}.$ Consequently, one gets 
$$\frac{\|P_k^{(\alpha)}\|_{\infty}}{\sqrt{h_k^{(\alpha)}}} \le 
\frac{\pi^{\frac{1}{4}} k^{\alpha+\frac{1}{2}}e^{\alpha-\frac{1}{4}}}{\sqrt{2}}
\left(1+\frac{1}{\alpha+\frac{1}{2}} \right)^{\alpha+\frac{1}{2}}
\leq \frac{\pi^{\frac{1}{4}}e^{\alpha+\frac{3}{4}}}{\sqrt{2}} k^{\alpha+\frac{1}{2}}.  $$
Finally, for  $k=0$,  we have
 $\|\widetilde{P}_0^{(\alpha)}\|_{\infty}^2
=\frac{1}{h_0^{(\alpha)}} $ which is, according to Lemma~\ref{lm:h}, upper bounded by $C^2(\alpha)$.

\begin{corollary}
	Under the same hypothesis and notations of the previous proposition, for any $\alpha \geq -\frac{1}{2}$ and for any 
	integer $N\geq 1,$ we have
\begin{equation}\label{bound_Sum}
\sum_{k=0}^{N} \|\widetilde{P}_k^{(\alpha)}\|_{\infty} \le
	\eta_{\alpha} 
	\frac{\left(N+1\right)^{\alpha+\frac{3}{2}}}{\alpha+\frac{3}{2}},\qquad \eta_\alpha=
	 \frac{\pi^{\frac{1}{4}} e^{\alpha+\frac{3}{4}}}{\sqrt{2}}.
\end{equation}
	\end{corollary}
\proof~ From the previous proposition, we can write
\begin{equation}
\begin{aligned}
  \sum_{k=0}^{N} \|\widetilde{P}_k^{(\alpha)}\|_{\infty} & =
     \|\widetilde{P}_0^{(\alpha)}\|_{\infty} +\sum_{k=1}^{N} \|\widetilde{P}_k^{(\alpha)}\|_{\infty} 
      \le
     \frac{\eta_{\alpha}}{\alpha+\frac{3}{2}}
     +
     \sum_{k=1}^N\eta_{\alpha} k^{\alpha+\frac{1}{2}} \\
     &\le  \eta_{\alpha} \int_0^{N+1} x^{\alpha+\frac{1}{2}}dx=\eta_{\alpha}\frac{\left(N+1\right)^{\alpha+\frac{3}{2}}}{\alpha+\frac{3}{2}}\,.
\end{aligned}
\end{equation}
\qed
~\\

Let 
$L^2(I,w_{\alpha})$ be the  Hilbert space associated to the inner  product $<f,g>_\alpha=\int_I f \cdot g \cdot w_{\alpha}(x)\, dx.$ Note that the family $\left\{\widetilde{P}_k^{(\alpha)}, k \ge 0 \right\}$ is an orthonormal basis of $L^2(I,w_{\alpha}).$

\section{An ANOVA type space based on multivariate Jacobi polynoimals}

In this paragraph, we describe a reduced size multidimensional polynomials space. The construction of this space is based on combining the ANOVA decomposition technique see for example \cite{Huang, Lin2, Potts, Saltelli, Sobol, Stone2}  and the total degree polynomial space, see for example \cite{Blatman}. For this purpose, let $\mathcal D =[[1,d]]=\{1,2,\ldots,d\}$ and we will adopt  the notations $\pmb u \subset \mathcal D$ for subsets of $\mathcal D,$ $x_{\pmb u}=(x_i)_{i\in \pmb u}$ and $|\pmb u|$ for the length of the vector $\pmb u$. 
For a given $\pu \subset \mathcal{D}$, we let $F_{\pu}$ denote the  subset of $\Z^d$ defined by
$$ F_{\pu}=\left\{ \pk \in \Z^d / k_{\pu^c}=0_{\Z^{d-|u|}} \, ;
\ k_j\neq 0 \ \forall j\,  \in \pu 
\right\}\,.$$

First, we describe a Jacobi polynomials  orthonormal  basis of  $L^2(I^d,\pw_{\alpha})$ based on the ANOVA decomposition. Here, the $d-$variate weight function $\pw$ is defined on 
$I^d$ by  
$$\forall \ \px=(x_1,\ldots,x_d) \in I^d, \ \pw_{\alpha}(\px)=\prod_{i=1}^{d} w_{\alpha}(x_i)=\left(\prod_{i=1}^d \left(1-x_i^2\right)\right)^{\alpha}\,.$$
The usual  inner product associated with  $L^2(I^d,\pw_{\alpha})$ is defined by 
$$<f,g>_{\pw}=\int_{[-1,1]^d} f(\pmb x)  \cdot g(\pmb x) \,  \pw_{\alpha}(\pmb x) \, d\pmb x.$$
 Let $\Phi_{ \pmb{0}}^{(\alpha)}(\pmb{x})=\frac{1}{\sqrt{h_0^{(\alpha)}}^{d}}$ and  
 for $\pmb{u} \subset \mathcal{D}$ such that $|\pu| \ge 1$ and $\pmb{k} \in F_{\pu}$, let
\begin{equation}
\Phi_{\pmb{u,k}}^{(\alpha)}(\pmb{x}):=
 \frac{1}{\sqrt{h_0^{(\alpha)}}^{d-|\pu|}}
 \prod_{i \in \pmb{u} } \widetilde{P}_{k_i}^{(\alpha)}(x_i) \,. 
\end{equation} 
 
 \begin{lemma}
 \label{lm:orthonormality}
     The family
$$
\left\{ \Phi_{\pmb{0}}^{(\alpha)} \right\} \cup
\left\{\Phi_{\pmb{u,k}}^{(\alpha)} , \  \pmb{u} \subset \mathcal{D}, \  \ \pmb{k} \in F_{\pu} \right\}$$
is  an orthonormal basis  of $L^2(I^d,\pw_{\alpha})$. 
 \end{lemma}
 
 \proof~ Since $$\mbox{Span}\left\{ \Phi_{\pmb{0}}^{(\alpha)} \right\} \cup
\left\{\Phi_{\pmb{u,k}}^{(\alpha)} , \  \pmb{u} \subset \mathcal{D}, \  \ \pmb{k} \in F_{\pu} \right\}=\mbox{Span}\left\{\prod_{i=1}^d P^{(\alpha)}_{k_i},\, k_i \geq 0\right\}$$ and since this later is dense in $L^2(I^d,\pw_\alpha),$ then it suffices to 
 establish the orthornormality of the vectors
 $\Phi_{\pmb{0}}^{(\alpha)}$ and $\Phi_{\pmb{u,k}}^{(\alpha)}$, we consider the following three cases.
 \begin{enumerate}
     \item {Computation of $<\Phi_{\pmb{u,k}}^{(\alpha)},\Phi_{\pmb{u,k}}^{(\alpha)}>_{\alpha}$}. 
     We first assume that   $|\pu|\ge 1$, then we have
     \begin{equation}
         \begin{aligned}
          <\Phi_{\pmb{u,k}}^{(\alpha)},\Phi_{\pmb{u,k}}^{(\alpha)}>_{\alpha} & = \int_{[-1,1]^d} \Phi_{\pmb{u,k}}^{(\alpha)}(\px) \prod_{i=1}^d w_{\alpha}(x_i)dx_1\ldots dx_d \\
        & =  \int_{[-1,1]^d} \left(\widetilde{P}_0^{(\alpha)}\right)^{2d-2|\pu|}
          \prod_{i \in \pu} \left(\widetilde{P}_{k_i}^{(\alpha)}(x_i)\right)^2
          \prod_{i=1}^d w_{\alpha}(x_i)dx_1\ldots dx_d \\
          & = \left(\prod_{i \in \pu} \int_{-1}^1 
          \left(\widetilde{P}_{k_i}^{(\alpha)}(x_i)\right)^2 w_{\alpha}(x_i) dx_i\right)
          \left( \prod_{1\le j \le d ; j \notin \pu}
          \int_{-1}^1 
          \left(\widetilde{P}_{0}^{(\alpha)}(x_j)\right)^2 w_{\alpha}(x_j) dx_j          \right)=1\,.
         \end{aligned}
     \end{equation}
     In a similar manner, we have $<\Phi_{\pmb{0}}^{(\alpha)},\Phi_{\pmb{0}}^{(\alpha)}>_{\alpha}=1$. 

     \item Computation of $<\Phi_{\pmb{u,k}}^{(\alpha)},\Phi_{\pmb{u,l}}^{(\alpha)}>_{\alpha}$ in the case where $|\pu| \ge 1$ and $\pk \neq \pmb{l}$ \\
     As $\pk \neq \pmb{l}$ then $\exists \ i_0 \in \pu$ such that $\pk_{i_0} \neq \pmb{l}_{i_0} $\,.\\
     \begin{equation}
     \label{eq:orth2}
     \begin{aligned}
      <\Phi_{\pmb{u,k}}^{(\alpha)},\Phi_{\pmb{u,l}}^{(\alpha)}>_{\alpha} & =\int_{[-1,1]^d} \Phi_{\pmb{u,k}}^{(\alpha)}(\px)\Phi_{\pmb{u,l}}^{(\alpha)}(\px) \pmb{w}_{\alpha}(\px)d\px\\
     & =\left(\frac{1}{h_0^{\alpha}}\right)^{d-|\pu|}
     \left(\int_{-1}^1 \widetilde{P}_{k_{i_0}}^{(\alpha)}(x_{i_0})
     \widetilde{P}_{l_{i_0}}^{(\alpha)}(x_{i_0})w_{\alpha}(x_{i_0})dx_{i_0}\right)
     \prod_{j \neq i_0 ; j \in \pu} 
     \left(\int_{-1}^1 \widetilde{P}_{k_{j}}^{(\alpha)}(x_{j})
     \widetilde{P}_{l_{j}}^{(\alpha)}(x_{j})w_{\alpha}(x_{j})dx_{j}\right)\\
     &=0
     \end{aligned}
     \end{equation}
     \item Computation of $<\Phi_{\pmb{u,k}}^{(\alpha)},\Phi_{\pmb{v,l}}^{(\alpha)}>_{\alpha}$.     Let us consider the case where $\pu \neq \pmb{v}$, $|\pu| \ge 1$ and $|\pmb{v}|\ge 1$\\
     As $\pu \neq \pmb{v}$ then : (1) $\pk \neq \pmb{l}$ and 
     (2) $\exists \ i_0 \in \pu \backslash \pmb{v}$ (or $\pmb{v} \backslash v$). Without loss of generality, we will suppose that $ i_0 \in \pu \backslash \pmb{v}$.
     \begin{equation}
     \label{eq:orth3}
     \begin{aligned}
      & <\Phi_{\pmb{u,k}}^{(\alpha)},\Phi_{\pmb{v,l}}^{(\alpha)}>_{\alpha} \\
      &=\int_{[-1,1]^d} \Phi_{\pmb{u,k}}^{(\alpha)}(\px)\Phi_{\pmb{v,l}}^{(\alpha)}(\px) \pmb{w}_{\alpha}(\px)d\px\\
     & = \int_{[-1,1]^d} \left( \prod_{i \in \pu} \widetilde{P}_{k_{i}}^{(\alpha)}(x_{i}) \right)
     \left(\frac{1}{h_0^{\alpha}}\right)^{\frac{d-|\pu|}{2}}
    \left( \prod_{j \in \pmb{v}} \widetilde{P}_{l_{j}}^{(\alpha)}(x_{j})
     \left(\frac{1}{h_0^{\alpha}}\right)^{\frac{d-|\pmb{v}|}{2}}\right)
     \prod_{s=1}^d w_{\alpha}(x_s)dx_s\\
     &=\left(\int_{-1}^1 \widetilde{P}_{k_{i_0}}^{(\alpha)}(x_{i_0}) \frac{w_{\alpha}(x_{i_0})}{\sqrt{h_0^{(\alpha)}}}dx_{i_0}\right)
     \int_{[-1,1]^{d-1}} 
     \frac{\left( \prod_{i \in \pu \backslash i_0; j \in \pmb{v}} \widetilde{P}_{k_{i}}^{(\alpha)}(x_{i}) 
 \widetilde{P}_{l_{j}}^{(\alpha)}(x_{j}) \right)}
 {\sqrt{h_0^{(\alpha)}}^{2d-|\pu|-|\pmb{v}|-1}}
       \prod_{s \in [[1,d]]\backslash i_0} w_{\alpha}(x_s) dx_s=0\\
     \end{aligned}
     \end{equation}
       Similarly to (\ref{eq:orth3}), we have
     in the case where $|\pmb{v}|=0$
     $<\Phi_{\pmb{u,k}}^{(\alpha)},\Phi_{\pmb{0}}^{(\alpha)}>_{\alpha}=0$.
 \end{enumerate}
 \qed

Next, we consider the following reduced size polynomial space,
defined for $m \in [[1,N]]$  on which we will construct our estimator. 
\begin{equation}
	\label{eq:Pspace}
\mathcal{P}_{N,m,d}= Span\left\{ \prod_{i \in \pmb u} \widetilde{P_{k_i}}^{(\alpha)}(x_i) \ ; \ \pu \subset \mathcal{D} \ ;
\ \ |\pu| \le m ; \ \pk \in F_{\pu}, \ \ ||\pk||_1 \le N  \right\}\,.
\end{equation}
We will call this space  ANOVA type polynomial space.
The dimension of the space $\mathcal{P}_{N,m,d}$  as well as an estimate of its upper bound are  given by  the following proposition.
\begin{proposition}\label{lm:DimOrder}
For any positive integers $1\leq m\leq \min(d, N),$ we have 
	\begin{equation}
	\label{eq1:dim}
	    M_{N,m,d}=\dim \mathcal{P}_{N,m,d}= \sum_{k=0}^m {d \choose k}{N \choose k}.  
	\end{equation}
Moreover, we have 
\begin{equation}
	\label{eq2:dim}
M_{N,m,d} \leq \left\{\begin{array}{ll} \frac{\pi}{ e}  \Big(1+\frac{m}{\min(d,N)-m+\frac{1}{2}}\Big)^{2\big(\max(d,N)-m+\frac{1}{2}\big)}
\cdot \left( \frac{\max(d,N)+\frac{1}{2}}{m+\frac{1}{2}} \right)^{2m}
&
\mbox{ if } m\leq \frac{1}{2} \min(d,N),\\
& \\
2^{\min(d,N)} \Big(1+\frac{m}{\max(d,N)-m}\Big)^{\max(d,N)-m}\Big(\frac{\max(d,N)}{m}\Big)^m&\mbox{ if } \frac{1}{2} \min(d,N) \leq m< \min(d,N)
\end{array}\right.
\end{equation}
and, if $m=\min(d,N)$,
\begin{equation}
    M_{N,m,d}={N+d \choose m} \le \sqrt{\frac{\pi}{(2d+2N+1) e}} \Big(1+\frac{N}{d+\frac{1}{2}}\Big)
    ^{d+\frac{1}{2}}\Big(1+\frac{d}{N+\frac{1}{2}}\Big)^{N+\frac{1}{2}} 
\end{equation}
\end{proposition}
	
\proof~ We first check \eqref{eq1:dim}. 
For this purpose, we first consider the two special cases $m=2$ and $m=3,$. Then, we check the previous identity for any $m\leq \min(d,N).$  For $m=2,$ the space $\mathcal{P}_{N,2,d}$ can be rewritten as $$\mathcal{P}_{N,2,d}=H_1^N \oplus H_2^N\,, $$
with 
$$H_1^N=Span\{ \widetilde{P}_{l}^{(\alpha)}(x_i) , \ 0 \le  l \le N , \ \ 1\le i \le d  \}\,$$
and 
$$H_2^N=Span\{ \widetilde{P}_{k_1}^{(\alpha)}(x_i)\widetilde{P}_{k_2}^{(\alpha)}(x_j) , \ (k_1,k_2) \in [[1,N-1]]^2 , \ 
 2 \le  k_1+k_2 \le N , \ \ 1\le i < j \le d  \}\,.$$
It is not difficult to check that the dimensions of $H_1^N$ 
and  $H_2^N$ are given by 
$$ \dim H_1^N = 1+ N d = \sum_{k=0}^1 {d \choose k}{N \choose k}, \qquad \dim H_2^N =\frac{d(d-1)N(N-1)}{4} = {d \choose 2}{N \choose 2}.$$
Hence, 
$$\dim \mathcal{P}_{N,2,d} =M_{N,2,d}=\sum_{k=0}^2 {d \choose k}{N \choose k}.$$
In a similar manner; for $m=3,$ we have $\mathcal{P}_{N,3,d}=\mathcal{P}_{N,2,d} \oplus H_3^N\,, $ with 
$$H_3^N=Span\Big\{\prod_{i=1}^3 \widetilde{P}_{k_i}^{(\alpha)}(x_{j_i}), \
(k_1,k_2,k_3) \in [[1,N-2]]^3, \ 
 3 \le  k_1+k_2+k_3 \le N , \ \ 1\le j_1 < j_2 < j_3\le d  \Big\}\,.$$
 Note that there exist ${d\choose 3}$ different  $3-$tuples 
 $(x_{j_1},x_{j_2},x_{j_3})$ in $H_3^N.$ Moreover, for each
 such a tuple, there correspond ${\displaystyle \frac{1}{2} \sum_{k=2}^{N-1} k(k-1)= {N\choose 3}}$ different polynomials products 
 ${\displaystyle \prod_{i=1}^3 \widetilde{P}_{k_i}, \, 3\leq \sum_{i=1}^3 k_i\leq N.}$ Consequently, we have 
 ${\displaystyle \dim H_3^N = {d\choose 3}{N\choose 3}}$ and 
 $$\dim \mathcal{P}_{N,3,d}=\dim \mathcal{P}_{N,2,d}+\dim H_3^N =\sum_{k=0}^3 {d \choose k}{N \choose k}.$$
 Continuing in this manner, one gets the recurrence formula 
 $$\dim \mathcal{P}_{N,m,d}=\dim \mathcal{P}_{N,m-1,d}+\dim H_m^N =\sum_{k=0}^{m-1} {d \choose k}{N \choose k}+{d \choose m}{N \choose m}=\sum_{k=0}^{m} {d \choose k}{N \choose k}.$$
Next, to prove \eqref{eq2:dim}, we proceed as follows.
 By using \eqref{eq:Gamma}, one gets  for any integers $n\geq k\geq 0$
\begin{eqnarray}\label{eq:combinaison}
{n \choose k} &=& \frac{\Gamma(n+1)}{\Gamma(k+1) \Gamma(n-k+1)}\leq \sqrt{\frac{\pi}{2 e^2}}
\frac{\left(\frac{n+\frac{1}{2}}{e}\right)
^{n+ \frac{1}{2} }}
{\left(\frac{k+\frac{1}{2}}{e}\right)
^{k+ \frac{1}{2} }
\left(\frac{n-k+  \frac{1}{2}}{e}\right)
^{n-k+ \frac{1}{2}  }} \nonumber\\
&\leq& \sqrt{\frac{\pi}{2 e}}\frac{1}{
\sqrt{k+\frac{1}{2}}}
\left(1+\frac{k}{n-k+\frac{1}{2}}\right)^{n-k+\frac{1}{2}}
\left(\frac{n+\frac{1}{2}}{k+\frac{1}{2}}\right)^{k}\,.
\end{eqnarray}
Since for $1\leq m\leq \frac{1}{2} \min(d,N),$ the finite sequences $\Big\{{d\choose i}\Big\}_{0\leq i\leq m}, \Big\{{N\choose i}\Big\}_{0\leq i\leq m}$ are increasing, then by using \eqref{eq:combinaison}, one gets
\begin{eqnarray*}
M_{N,m,d} &\leq & (m+1) {d\choose m} {N\choose m} = (m+1)\, \frac{ \Gamma(d+1) \cdot \Gamma(N+1)}{\Gamma(m+1) 
\Gamma\left(d-m+\frac{1}{2}\right)\cdot \Gamma(m+1) \Gamma\left(N-m+\frac{1}{2}\right)}\\
&\leq & \frac{\pi}{ e}  \Big(1+\frac{m}{\min(d,N)-m+\frac{1}{2}}\Big)^{2\big(\max(d,N)-m+\frac{1}{2}\big)}
\cdot \left( \frac{\max(d,N)+\frac{1}{2}}{m+\frac{1}{2}} \right)^{2m}\,.
\end{eqnarray*}
In a similar manner, for $ \frac{1}{2} \min(d,N) \leq m <  \min(d,N),$ we have 
\begin{eqnarray*}
M_{N,m,d} &= & \sum_{k=0}^{m} {d\choose k} {N\choose k} \leq 
\sum_{k=0}^{\min(d,N)} {\min(d,N)\choose k} \cdot \sum_{k=0}^{m}{\max(d,N)\choose k}\\
&\leq& 2^{\min(d,N)} \frac{\max(d,N)^{\max(d,N)}}{m^m ((\max(d,N)-m)^{\max(d,N)-m}}.
\end{eqnarray*}
This last inequality is a consequence of the inequality,  $$\sum_{k=0}^m {n\choose k} \leq \frac{n^n}{m^m (n-m)^{n-m}},\qquad \forall\,\, 0\leq m< n.$$
The previous inequality is a direct consequence of the following upper bound for the binomial coefficient, see [\cite{Cover}, p.353] 
$$ {n\choose k} \leq {\displaystyle 2^{n H\big(k/nbig)}},$$
where for $0<p<1,$ $H(p)$ is the binary entropy function, given by 
${\displaystyle H(p)=- p \log_2(p) -(1-p) \log_2(1-p).}$
For the last case where $m=\min(d,N)$, using 
Chu-Vandermonde's identity, see for example \cite{Roy}, we get
$$ \sum_{k=0}^{m} {d\choose k} {N\choose k} =\sum_{k=0}^{\min(d,N)} { \min(d,N) \choose k} {\max(d,N)\choose k} 
=\sum_{k=0}^{\min(d,N)} { \min(d,N) \choose \min(d,N)- k} {\max(d,N)\choose k} = 
{N+d \choose m}.$$
Using \ref{eq:Gamma}, we get
\begin{equation}
     {N+d \choose m}=\frac{\Gamma(N+d+1)}{\Gamma(d+1)\Gamma(N+1)}  \le \sqrt{\frac{\pi}{2e}}\frac{1}{\sqrt{N+d+\frac{1}{2}}} \left( \frac{N+d+\frac{1}{2}}{d+\frac{1}{2}}\right)^{d+\frac{1}{2}}
     \left( \frac{N+d+\frac{1}{2}}{N+\frac{1}{2}}\right)^{N+\frac{1}{2}}\,.
\end{equation}
\begin{remark}
 For the particular  case $m=d$, the space 
 $\mathcal P_{N,d,d} $ is the usual total degree polynomial of order $N.$ 
 From Proposition~\ref{lm:DimOrder}, we have $M_{N,m,d}  = {N+d \choose d}$ which recovers the well known  dimension of the total degree polynomial space. 
 \end{remark}

\section{Stability of the NP estimator}
\label{Sec:stability}
In this section, we first describe our proposed Least squares NP regression estimator based on the use of the  orthonormal $d-$variate polynomials $\Phi_{\pmb{u,k}}^{(\alpha)}.$ Then, we prove the stability of the proposed estimator under the condition that the random  sampling covariates follow a $d-$dimensional Beta distribution. \\

Recall that the regression system at hand is given by 
\begin{equation}
	\label{eq:Pb}
	Y_i=f(X_i)+\varepsilon_i, \  i=1,\ldots,n \,.
\end{equation}
Here, $f$ is the regression function to be approximated, the $X_i \in \mathbb R^d$ are the covariates which we assume to be i.i.d. random   vectors
and the 
 $(\varepsilon_i)_{1\le i \le n}$ are centered i.i.d. random variables with a finite variance 
$\sigma^2=\mathbb E[\varepsilon_1^2]$. Also we assume that the $X_i$ are independent from the $\varepsilon_j.$ In the sequel, we assume that for a real $\alpha\geq -\frac{1}{2},$ the set $\{X_i,\, i=1,\ldots,n\}$ is a  random sampling set with the $X_i$ following a $d-$variates $\mbox{Beta}(\alpha+1,\alpha+1)-$distribution.  
For two positive integers $m\leq d$ and $N$,
we  build an estimator $\widehat{f}_{N,n,m}^{(\alpha)}$ of the regression function $f$ which is obtained using the approximation of $f$ by its projection over $\mathcal{P}_{N,m,d}$. For $\pu  \subset \mathcal{D}$ with $|u|\ge 1$ and $\pk \in F_{\pu}$, let $\Phi_{k_{\pu}}^{(\alpha)}:=\Phi_{\pmb{u,k}}^{(\alpha)}$.
Let $K_{N,0}:=\left\{ (0,\ldots,0) \in \R^d \right\}$  and 
for $p \in [[1,m]]$, let $K_{N,p}$ be the subsets defined as follows
$$K_{N,p}:=
\left\{ k_{\pu} \in   F_{\pu} : 
|\pmb{u}|=p  \ ; \  ||k_{\pu}||_1 \le N
\right\}\,\qquad K_{N,m,d}:=  \bigcup_{p=0}^m K_{N,p}. $$
Let $g : K_{N,m,d} \mapsto [[1,M_{N,m,d}]]$ be a correspondence (order) defined on the indexes $(\pk,\pu)$ of our basis : $g(\pk,\pu) \in [[1,M_{N,m,d}]]$. Then, we can introduce the notation $$ \Psi_j^{(\alpha)}:=\Phi_{g^{-1}(j)}^{(\alpha)}, \ \mbox{with} \  i \in [[1,M_{N,m,d}]]\,.$$
Using these definitions, our NP regression estimator is given by 
\begin{equation}\label{estimator}
	\widehat{f}_{N,n,m}^{(\alpha)}(x)= \sum_{k_{\pu} \in  K_{N,m,d}}
	\widehat{C}_{k_{\pu}}(f) \Phi_{k_{\pu}}(x)
	= \sum_{j=1}^{M_{N,m,d}} \widehat{C}_{j}(f) \Psi_j(x),\quad x\in I^d.
\end{equation}
Assuming that $\widehat{f}_{N,n,m}^{(\alpha)}$ satisfies  $\widehat{f}_{N,n,m}^{(\alpha)}(X_i)=Y_i$
and multiplying the previous equation by $\frac{\left[2^{2\alpha+1}\B(\alpha+1,\alpha+1)\right]^\frac{d}{2}}{n^{\frac{1}{2}}}$, we obtain
\begin{equation}
	\sum_{j=1}^{M_{N,m,d}} \widehat{C}_{j}(f) 	\frac{\left[2^{2\alpha+1}\B(\alpha+1,\alpha+1)\right]^\frac{d}{2}}{n^{\frac{1}{2}}}\Psi_{j}(X_i)=
	\frac{\left[2^{2\alpha+1}\B(\alpha+1,\alpha+1)\right]^\frac{d}{2}}{n^{\frac{1}{2}}}
	Y_i.
\end{equation}
This system can be rewritten  as a system of linear equations where the unknown is the expansion coefficients vector $\widehat{C}_{N,n,m}^{(\alpha)}=
\left(\widehat{C}_{j}\right)_{1 \le j \le M_{N,m,d}}^T$:
$$ F_{N,n,m}^{(\alpha)}.\widehat{C}_{N,n,m}^{(\alpha)}=
\frac{\left[2^{2\alpha+1}\B(\alpha+1,\alpha+1)\right]
^\frac{d}{2}}{n^{\frac{1}{2}}}
Y.$$
Consider the positive definite random matrix (to be checked later on),
\begin{equation}\label{matrixG}
G_{N,n,m}^{(\alpha)}=
\left(F_{N,n,m}^{(\alpha)}\right)^T \cdot F_{N,n,m}^{(\alpha)}.
\end{equation}
Then, we have
\begin{equation}
\label{eq:C_hat}
\widehat{C}_{N,n,m}^{(\alpha)}=
\left(G_{N,n,m}^{(\alpha)}\right)^{-1}
\left(F_{N,n,m}^{(\alpha)}\right)^T.\frac{\left[2^{2\alpha+1}\B(\alpha+1,\alpha+1)\right]^\frac{d}{2}}{n^{\frac{1}{2}}} Y,    
\end{equation}
with
\begin{equation}\label{eq:randomamatrix}
F_{N,n,m}^{(\alpha)}(i,j)=\frac{\left[2^{2\alpha+1}\B(\alpha+1,\alpha+1)\right]^\frac{d}{2}}{n^{\frac{1}{2}}}.\Psi_{j}^{(\alpha)}(X_i).
\end{equation}

Next, we show that our proposed NP regression estimator is stable in the sense that the  random matrix $G_{N,n,m}^{(\alpha)}$ is well conditioned with respect to the $2-$norm. As the random matrix 
$G_{N,n,m}^{(\alpha)}$ is positive definite, its
condition number denoted by $\kappa_2\left(G_{N,n,m}^{(\alpha)}\right)$  is given by 
\begin{equation}\label{condition_number}
\kappa_2\left(G_{N,n,m}^{(\alpha)}\right)=
\frac{\lambda_{\max}\left(G_{N,n,m}^{(\alpha)}\right)}
{\lambda_{\min}\left(G_{N,n,m}^{(\alpha)}\right)}.
\end{equation}

Recall that  the i.i.d. random samples 
$(X_i)_{1 \le i \le  n}$ follow the $d-$variate Beta distribution on $I^d$
with the density function $\pmb{g}_{\alpha} $ given by 
\begin{equation} \label{eq:density}
\pmb{g}_{\alpha}(\px)=\frac{1}{\left[2^{2\alpha+1}\B(\alpha+1,\alpha+1)\right]^d}\pw_{\alpha}(x){\mathds 1}_{I^d}(\px)=\frac{1}{\left(h_0^{(\alpha)}\right)^d}\pw_{\alpha}(\px){\mathds 1}_{I^d}(\px)\,. \end{equation}
Before stating the theorem about the condition number upper bound (Theorem~\ref{th:CondNumber}), we need  the following technical lemma.
\begin{lemma}
	\label{lm:Bounds}
	For $\alpha \ge -\frac{1}{2}$, $N \in \N$ and a positive  integer $m\leq d,$ let
	\begin{equation}\label{Lm1:Bounds}
	    D(N,\alpha):=
	    \begin{cases}
	  \pi^{\frac{3}{4}}  N^{(\alpha+\frac{1}{2})}
e^{(\alpha+\frac{3}{4})} \ \ 
 \mbox{if} \ \  \alpha >-\frac{1}{2}\\
	    2 \ \ \mbox{if} \ \ \alpha=-\frac{1}{2}\,.
	    \end{cases}
	\end{equation}
	Then, we have 
	\begin{equation}
	\label{eq:boundphi0}
				  \left\| \Phi_{\pmb{0}}^{(\alpha)}\right\|_{\infty}
		 = \left(\frac{1}{h_0^{(\alpha)}}\right)^{\frac{d}{2}}
		 \,.
	\end{equation}
Moreover,   for  $1 \le |\pu|\le m$ and $||\pk_{\pu}||_1 \le N$, we have  
\begin{equation}
\label{eq:boundphik}
	\left\| \Phi_{k_{\pu}}^{(\alpha)}\right\|_{\infty}\le 
\left(\frac{1}{h_0^{(\alpha)}}\right)^{\frac{d}{2}} \cdot [D(N,\alpha)]^m\,.
\end{equation}
\end{lemma}
\proof~  
Equality (\ref{eq:boundphi0})
is  immediate since $\Phi_0^{(\alpha)}=\left[\widetilde{P}_0^{(\alpha)}\right]^d=
\left(\frac{1}{h_0^{(\alpha)}}\right)^{\frac{d}{2}}
$. Also from the definition of $\Phi_{k_{\pu}}^{(\alpha)}$, we have
$$\|\Phi_{k_{\pu}}^{(\alpha)}\|_{\infty} 
\le \|\widetilde{P}_0^{(\alpha)} \|_{\infty}^{d-|\pu|} 
\prod_{i \in \pu} \|  \widetilde{P}_{k_i}^{(\alpha)}\|_{\infty}
=
\left(\frac{1}{h_0^{(\alpha)}}\right)^
{\frac{d-|\pu|}{2}}
\prod_{i \in \pu} \|  \widetilde{P}_{k_i}^{(\alpha)}\|_{\infty}
\,.
$$
Next, for the case where $\alpha > -\frac{1}{2}$ and by using 
 Proposition \ref{prop:boundsalpha}
  and Lemma~\ref{lm:h}, one gets :
\begin{equation}
\begin{aligned}
 \|\Phi_{k_{\pu}}^{(\alpha)}\|_{\infty} 
& \le  
\left(\frac{1}{h_0^{(\alpha)}}\right)^
{\frac{d-|\pu|}{2}}
\prod_{i \in \pu}
\left[ \frac{\pi^{\frac{1}{4}} N^{\alpha+\frac{1}{2}}e^{\alpha+\frac{3}{4}}}{\sqrt{2}}
\right] \\
~ & =
 \left(\frac{1}{h_0^{(\alpha)}}\right)^
{\frac{d}{2}}
\left(h_0^{(\alpha)}\right)^
{\frac{|\pu|}{2}}
 \left[ \frac{\pi^{\frac{|\pu|}{4}} N^{(\alpha+\frac{1}{2})|\pu|}e^{(\alpha+\frac{3}{4})
 |\pu|}}{\sqrt{2}^{|\pu|}}
\right]\\
~ & \le  
 \left(\frac{1}{h_0^{(\alpha)}}\right)^
{\frac{d}{2}}
\left[ 
\left(\pi^3\right)^
{\frac{|\pu|}{4}}  N^{(\alpha+\frac{1}{2})|\pu|}
e^{(\alpha+\frac{3}{4})|\pu|}
\right]\,. \\
~ & \le  
 \left(\frac{1}{h_0^{(\alpha)}}\right)^
{\frac{d}{2}}
\left[ 
\pi^{\frac{3}{4}}  N^{(\alpha+\frac{1}{2})}
e^{(\alpha+\frac{3}{4})}
\right]^m\,. \\
\end{aligned}
\end{equation}
Finally, for the case   $\alpha=-\frac{1}{2}$ and by 
using again Proposition~\ref{prop:boundsalpha}
and Lemma~\ref{lm:h}, one gets
\begin{equation}
\begin{aligned}
 \|\Phi_{k_{\pu}}^{(\alpha)}\|_{\infty} 
& \le   \left(\frac{1}{h_0^{(\alpha)}}\right)^
{\frac{d}{2}}
\left(h_0^{(-\frac{1}{2})}\right)^
{\frac{|\pu|}{2}}
\prod_{i \in \pu}
\left[ \frac{2}{\sqrt{\pi}}
\right] \\
~ & \le
 \left(\frac{1}{h_0^{(\alpha)}}\right)^
{\frac{d}{2}}
\pi^
{\frac{|\pu|}{2}}
\left[ \frac{2}{\sqrt{\pi}}
\right]^{|\pu|} \le  
 \left(\frac{1}{h_0^{(\alpha)}}\right)^
{\frac{d}{2}}
2^{m}
\,.
\end{aligned}
\end{equation}
\qed

In the following theorem, we show that with high probability, the $2-$norm condition number of the positive definite random matrix $G_{N,n,m}^{(\alpha)}$ is bounded by a convenient constant depending on a parameter $0<\delta <1.$
The results of this  theorem can be considered as a generalization of a similar result given in \cite{BDK} for tensor product multivariate Jacobi polynomials basis.

\begin{theorem}	\label{th:CondNumber}
Under the previous notation and hypotheses on $N,m,\alpha,$ let $0<\delta <1$ and let $G_{N,n,m}^{(\alpha)}$ be the random matrix given by \eqref{matrixG}. Then, we have 
	\begin{equation}
	\label{eq:CondUpBound}
	    \PP\left( \kappa_2\left(G_{N,n,m}^{(\alpha)}\right) \le \frac{1+\delta}{1-\delta}  \right)
	    \ge 1 -2M_{N,m,d}. \exp{\left(-\delta^2 \frac{n}{3
	    D^{2m}(N,\alpha)
	    M_{N,m,d}}\right)}.
	\end{equation}
	 Here,  the quantity $D(N,\alpha)$ is as defined by  Lemma~\ref{lm:Bounds} and $M_{N,m,d}$ is as given by Proposition 2.
\end{theorem}

\proof~From \eqref{eq:randomamatrix}, we have $F_{N,n,m}^{(\alpha)}F(p,i)=\frac{\left[2^{2\alpha+1}\B(\alpha+1,\alpha+1)\right]^\frac{d}{2}}{n^{1/2}}.\Psi_i^{(\alpha)}(X_p)$ and $G_{N,n,m}^{(\alpha)}
=\left(F_{N,n,m}^{(\alpha)}\right)^T \cdot F_{N,n,m}^{(\alpha)}\,,$ one gets 
$$ G_{N,n,m}^{(\alpha)}(i,j)=\sum_{p=1}^{n}  \frac{\left[2^{2\alpha+1}\B(\alpha+1,\alpha+1)\right]^d}{n}.\Psi_i^{(\alpha)}(X_p)\Psi_j^{(\alpha)}(X_p)=\sum_{p=1}^{n} 
\frac{\left( h_0^{(\alpha)}  \right)^d}{n}.\Psi_i^{(\alpha)}(X_p)\Psi_j^{(\alpha)}(X_p)
\,. $$
Consequently, one gets 
$$ \mathbb E[G_{N,n,m}^{(\alpha)}(i,j)]= \frac{\left( h_0^{(\alpha)}  \right)^d}{n} \sum_{p=1}^n
\mathbb E\left[\Psi_i^{(\alpha)}(X_p)\Psi_j^{(\alpha)}(X_p) \right]\,. $$
Now  to compute $\mathbb E\left[\Psi_i^{(\alpha)}(X_p)\Psi_j^{(\alpha)}(X_p) \right]$, we  let $i=g(\pk, \pu)$ and $j=g(\pmb{l},\pmb{v}),$ then by using Lemma~\ref{lm:orthonormality}, we obtain
\begin{equation}
    \begin{aligned}
    \mathbb  E\left[\Psi_i^{(\alpha)}(X_p)\Psi_j^{(\alpha)}(X_p) \right] & =
    \mathbb E\left[\Phi_{\pk, \pu}^{(\alpha)}(X_p)\Phi_{\pmb{l},\pmb{v}}^{(\alpha)}(X_p) \right]\\
     & = \int_{[-1,1]^d}
     \Phi_{\pk, \pu}^{(\alpha)}(\px)\Phi_{\pmb{l},\pmb{v}}^{(\alpha)}(\px)
     \frac{\prod_{s=1}^d w_{\alpha}(\px_s)d\px_s}{\left(h_0^{(\alpha)}\right)^d}\\
     & = \frac{1}{\left(h_0^{(\alpha)}\right)^d}
     <\Phi_{\pk, \pu}^{(\alpha)},\Phi_{\pmb{l},\pmb{v}}^{(\alpha)} >_{\alpha}=\frac{1}{\left(h_0^{(\alpha)}\right)^d}
     \delta_{(\pk, \pu),(\pmb{l},\pmb{v})}=\frac{1}{\left(h_0^{(\alpha)}\right)^d}
     \delta_{i,j}\,.
    \end{aligned}
\end{equation}
Hence, we have $\mathbb E\left[\Psi_i^{(\alpha)}(X_p)\Psi_j^{(\alpha)}(X_p) \right]=\frac{1}{\left(h_0^{(\alpha)}\right)^d}\delta_{i,j}$. Therefore $\mathbb E\big [G_{N,n,m}^{(\alpha)}\big ]$
is the identity matrix of dimension $M_{N,m,d}$.
Also, note  that
$$ G_{N,n,m}^{(\alpha)}=\sum_{p=1}^{n}H_p\,, $$
where $H_p(i,j)=\frac{\left(h_0^{(\alpha)}\right)^d}{n}.\Psi_i^{(\alpha)}(X_p)\Psi_j^{(\alpha)}(X_p)$.
We show that, for all $p \in [[1,n]]$, all eigenvalues of $H_p$ are non-negative. Moreover, we give an estimate for 
an  upper bound for the eigenvalues of  the  $H_p$.
The matrix $H_p$ can be written as $H_p=A^T_p.A_p$ where
 $A_p:=\frac{\left(h_0^{(\alpha)}\right)^{\frac{d}{2}}}{\sqrt{n}}.
 \left(	\Psi_1^{(\alpha)}(X_p), \ldots,\Psi_{M_{N,m,d}}^{(\alpha)}(X_p)\right)$.  The relation $H_p=A^T_p.A_p$
 implies that the eigenvalues of all the matrices $H_p$ are non-negative. Applying Gershgorin Theorem to $H_p$, one gets 
 $$ \lambda_{\max} (H_p) \le \max_{1 \le i \le M_{n,m,d}} \frac{\left(h_0^{(\alpha)}\right)^d}{n} 
 \left( \big |\Psi_i^{(\alpha)}(X_p)\big |^2 + \sum_{j \neq i} \big |\Psi_i^{(\alpha)}(X_p)\Psi_j^{(\alpha)}(X_p)\big |  \right).
 $$
Therefore
${\displaystyle \lambda_{\max}(H_p) \le \frac{ D(N,\alpha)^{2m}}{n} \cdot M_{N,m,d}}$ for all $l \in [[1,n]]$.
 Next, by applying the matrix Chernoff Theorem to  the matrix $G_{N,n,m}^{(\alpha)}$ as a sum of positive semi definite random matrices  $G_{N,n,m}^{(\alpha)}$ and  satisfying $0 \preceq H_p \preceq L \cdot  I_{M_{N,m,d}},$  $L:=\frac{\left(D(N,\alpha)\right)^{2m}\cdot M_{N,m,d}}{n}$ and as $\mu_{{\min}}=\lambda_{{\min}} \mathbb E\left(G_{N,n,m}^{(\alpha)}\right)=1$,
 $\mu_{{\max}}=\lambda_{{\max}} \mathbb E\left(G_{N,n,m}^{(\alpha)}\right)=1$, we get :
 \begin{equation}\label{eq:lambdamin}
     \PP
     \left(\lambda_{{\min}}\left(G_{N,n,m}^{(\alpha)}\right) \ge 1-\delta \right) 
     \ge 1- M_{N,m,d} \cdot \exp\left(-\frac{\delta^2n}{2
     [D(N,\alpha)]^{2m} \cdot 
     M_{N,m,d}}\right)
 \end{equation}
 and
 \begin{equation}\label{eq:lambdamax}
     \PP\left(\lambda_{{\max}}\left(G_{N,n,m}^{(\alpha)}\right)
     \le 1+\delta \right) 
     \ge 1- M_{N,m,d} \cdot \exp\left(-\frac{\delta^2n}{3
         [D(N,\alpha)]^{2m}
     \cdot
     M_{N,m,d}}\right)
 \end{equation}
 for $\delta \in (0,1)$. Finally, consider the events
 $A_1:= \left(\kappa_2\left(G_{N,n,m}^{(\alpha)}\right) \le \frac{1+\delta}{1-\delta}\right)$, $A_2=\left(\lambda_{min}\left(G_{N,n,m}^{(\alpha)}\right) \ge 1-\delta \right)$ and
 $A_3:=\left(\lambda_{\max}\left(G_{N,n,m}^{(\alpha)}\right) \le 1+\delta \right)$.
 Using the fact that $\PP(A_1^c) \le \PP(A_2^c) + \PP(A_3^c)$, we get
 $$\PP(A_1) \ge 
 1- 2 M_{n,m,d} \cdot \exp\left(-\frac{\delta^2n}{3
 [D(N,\alpha)]^{2m}. M_{N,m,d}}\right)
 \,.$$
\qed

\begin{remark} By using  \eqref{Lm1:Bounds} and \eqref{eq:CondUpBound}, one concludes that a convenient choice for the  value of the parameter $\alpha$ is given by $\alpha=-\frac{1}{2}.$ For this value, a given upper bound for $2-$norm condition number of the random matrix $G_{N,n,m}^{(\alpha)}$
is obtained with fewer number of $n,$ the number of random sampling points $X_i.$ This behaviour is illustrated by the numerical simulations given by Table 1 of the numerical examples section. 
\end{remark}

\section{$L^2-$ Risk error of the NP regression estimator}
\label{Sec:L2Risk}

In this section, we give an estimate for the $L^2-$risk error of our proposed NP regression estimator. For this purpose and it is done in \cite{Cohen},  we assume that there exists a constant $K_f$ such that
\begin{equation}
\label{eq:Hyp_f}
|f(\px)| \le K_f, \ \forall \ \px \ \in \ I^d\,.     
\end{equation}
We let $\widehat{F}_{N,n,m}^{(\alpha)}$ denote  the truncated version of the	estimator $\widehat{f}_{N,n,m}^{(\alpha)}$, given by 
\begin{equation}
\label{eq:F_Hat}
    \widehat{F}_{N,n,m}^{(\alpha)}(\px):=\mbox{Sign}\left( \widehat{f}_{N,n,m}^{(\alpha)}(\px)\right)\min 
    \left\{ K_f, |\widehat{f}_{N,n,m}^{(\alpha)}(\px)|\right\}\,.
\end{equation}
\begin{theorem} \label{th:l2risk}
Let $f$ be a function satisfying the  hypothesis (\ref{eq:Hyp_f}) and let
$\widehat{F}_{N,n,m}^{(\alpha)}$ be the truncated version of the	estimator $\widehat{f}_{N,n,m}^{(\alpha)}$, given by (\ref{eq:F_Hat}). For
 $\alpha \ge -\frac{1}{2}$ and $0 < \delta <1$, we have
\begin{multline}
    \mathbb E\left[ \| f- \widehat{F}_{N,n,m}^{(\alpha)}\|_{\alpha}^2\right] \le 
 4K_f^2
     \left(h_0^{(\alpha)}\right)^d
     M_{N,m,d} \cdot \exp\left(-\frac{\delta^2n}{2
     [D(N,\alpha)]^{2m} \cdot M_{N,m,d}}\right)\\
     +
 \frac{M_{N,m,d}}{n(1-\delta)^2}
     \bigg(
[D(N,\alpha)]^{2m} \|f-\Pi_{N,m}f\|_{\alpha}^2
+\sigma^2 \left(h_0^{(\alpha)}\right)^{d} \bigg) +
\|f-\Pi_{N,m}f \|_{\alpha}^2\,,
\end{multline}
 where the quantity  $ D(N,\alpha)$ is as defined by  Lemma~\ref{lm:Bounds}.
\end{theorem}
\proof~ As it is done in \cite{Cohen}, We write $\Omega$ as $\Omega_{+} \cup \Omega_{-},$ where 
$$ \Omega_{-}:=\{(X_1,\ldots,X_n) : \lambda_{{\min}}(G) < 1-\delta  \}\ \
\mbox{and} \ \ 
\Omega_{+}:=\{(X_1,\ldots,X_n) : \lambda_{{\min}}(G) \ge 1-\delta\}\,.
$$
Then 
$$ \mathbb E\left[ || f- \widehat{F}_{N,n,m}^{(\alpha)}||_{\alpha}^2\right]= \int_{\Omega_{-}} || f- \widehat{F}_{N,n,m}^{(\alpha)}||_{\alpha}^2  d\rho_n +
\int_{\Omega_{+}} || f- \widehat{F}_{N,n,m}^{(\alpha)}||_{\alpha}^2 d\rho_n\,,
$$
where $d\rho_n$ is the probability measure on $(I^d)^n$ given by the tensor product 
$$ d\rho_n(\pmb{u_1},\ldots,\pmb{u_n})=
\prod_{k=1}^n dg_{\alpha}(\pmb{u_k})\,,$$
with $g_{\alpha}$ given in (\ref{eq:density}).
To get an upper bound for $\int_{\Omega_{-}} \| f- \widehat{F}_{N,n,m}^{(\alpha)} \|_{\alpha}^2 d\rho_n$, we proceed as follows. 
From (\ref{eq:lambdamin}), we have
 $$ \PP\left(\lambda_{{\min}}
 \left(G_{N,n,m}^{(\alpha)}\right) < 1-\delta \right) 
     \le M_{N,m,d} \cdot\exp\left(-\frac{\delta^2n}{2
     [D(N,\alpha)]^{2m}\cdot
     M_{N,m,d}}\right)\,.$$
Then     
\begin{equation}
\begin{aligned}
 \int_{\Omega_{-}} \| f- \widehat{F}_{N,n,m}^{(\alpha)} \|_{\alpha}^2
 d\rho_n & \le 4K_f^2
     \int_{\Omega_{-}} d\rho_n =4K_f^2
     \left(h_0^{(\alpha)}\right)^d \PP\left(
     \lambda_{{\min}}
     \left(G_{N,n,m}^{(\alpha)}\right) < 1-\delta
     \right)\\
&      \le 
     4K_f^2
     \left(h_0^{(\alpha)}\right)^d
     M_{N,m,d} \cdot \exp\left(-\frac{\delta^2n}{2
     [D(N,\alpha)]^{2m} \cdot 
     M_{N,m,d}}\right)\,.
\end{aligned}
\end{equation}
For an estimate of an Upper bound of $\int_{\Omega_{+}} || f- \widehat{F}_{N,n,m}^{(\alpha)}||_{\alpha}^2 d\rho_n$, we use
 the fact that 
$|f-\widehat{F}_{N,n,m}^{(\alpha)}|$ is upper  bounded by $|f -\widehat{f}_{N,n,m}^{(\alpha)}|$ and the fact that  $(f-\Pi_{N,m} f)$ and  $(\Pi_{N,m} f-\widehat{f}_{N,n,m}^{(\alpha)})$
are orthogonal. Hence,  we get
\begin{equation}
    \begin{aligned}
     \int_{\Omega_{+}} \| f- \widehat{F}_{N,n,m}^{(\alpha)}\|_{\alpha}^2d\rho_n & \le \int_{\Omega_{+}} \| f- \widehat{f}_{N,n,m}^{(\alpha)}\|_{\alpha}^2 d\rho_n\\
     & \le \int_{\Omega_{+}} \| f- \Pi_{N,m}f\|_{\alpha}^2 d\rho_n
     + \int_{\Omega_{+}} \| \Pi_{N,m}f - \widehat{f}_{N,n,m}^{(\alpha)}||_{\alpha}^2 d\rho_n\,.
    \end{aligned}
\end{equation}
On the other hand,  let $C_{N,n,m}^{(\alpha)}$ the vector containing the coefficients of $\Pi_{N,m} f$ on the  basis
$\{\Psi_1^{(\alpha)},\ldots,\Psi_{M_{N,m,d}}^{(\alpha)}\}$,i.e.,
$\Pi_{N,m}f(X_i)=\sum_{j=1}^{M_{N,m,d} }  C_{j} \Psi_{j}^{(\alpha)}(X_i)$ for $i=1,\ldots,n$. This leads to 
$$C_{N,n,m}^{(\alpha)}=
\left(G_{n,N,m}^{(\alpha)}\right)^{-1}
\left(F_{N,n,m}^{(\alpha)}\right)^T.
\sqrt{\frac{\left(h_0^{(\alpha)}\right)^d}{n}}
\left(\Pi_{N,m}(f)(X_i)\right)_{1\le i \le n}\,.   
$$
The identity (\ref{eq:C_hat})  used for the computation of $\widehat{f}_{N,n,m}^{(\alpha)}$
 can be rewritten as:
$$\widehat{C}_{N,n,m}^{(\alpha)}=
\left(G_{N,n,m}^{(\alpha)}\right)^{-1}
\left(F_{N,n,m}^{(\alpha)}\right)^T \cdot \sqrt{\frac{\left(h_0^{(\alpha)}\right)^d}{n}} \left(f(X_i)+\epsilon_i)\right)_{1 \le i \le n}.$$
Combining the previous two equations, we get
$$ \widehat{C}_{N,n,m}^{(\alpha)}-C_{N,n,m}^{(\alpha)}=
\left(G_{N,n,m}^{(\alpha)}\right)^{-1}
(F_{N,n,m}^{(\alpha)})^T
\cdot \sqrt{\frac{\left(h_0^{(\alpha)}\right)^d}{n}} \cdot \left(g(X_i)+\epsilon_i\right)_{1\le i \le n}\,.
$$
where $g(X_i):=\left(f-\Pi_{N,m}(f)\right)(X_i)$.
From Parseval's equality, we have on $\Omega_{+}$
\begin{equation}
    \begin{aligned}
  \| \Pi_{N,m}(f) - \widehat{f}_{N,n,m}^{(\alpha)} \|_{\alpha}^2 & = 
  || \widehat{C}_{N,n,m}^{(\alpha)}-C_{N,n,m}^{(\alpha)}||_{l_2}^2    \\
&  \le 
\frac{\left(h_0^{(\alpha)}\right)^d}{n}
\left\|\left(G_{N,n,m}^{(\alpha)}\right)^{-1} \right\|_2^2 \cdot \left\| \left(F_{N,n,m}^{(\alpha)}\right)^T \cdot 
\left(g(X_i)+\epsilon_i\right)_{1\le i \le n} 
\right\|_{l_2}^2\\
& \le \frac{\left(h_0^{(\alpha)}\right)^d}{n(1-\delta)^2}
\left\| (F_{N,n,m}^{(\alpha)})^T\cdot 
\left(g(X_i)+\epsilon_i\right)_{1\le i \le n} 
\right\|_{l_2}^2\\
& = \frac{\left(h_0^{(\alpha)}\right)^{2d}}{n^2(1-\delta)^2}
\sum_{j=1}^{M_{N,m,d}} \sum_{k,l=1}^n \left(\Psi_j(X_k)(g(X_k)+\epsilon_k)\right)
\cdot
\left(\Psi_j(X_l)(g(X_l)+\epsilon_l)\right)\,.
    \end{aligned}
\end{equation}
Considering the expectations of both sides of the previous inequality and taking into account that\\
$<\Psi_j^{(\alpha)},g>_{\alpha}=0$  and that
the $\varepsilon_k$ are independent from the $X_k$ with  $\mathbb E(\varepsilon_k)=0$ and  $\mathbb E(\epsilon_k^2)=\sigma^2$, we get
\begin{equation}
    \begin{aligned}
   \mathbb  E\left[ \| \Pi_{N,m}(f) - \widehat{f}_{N,n,m}^{(\alpha)} \|_{\alpha}^2
    \right]
    & \le 
    \frac{\left(h_0^{(\alpha)}\right)^{2d}}{n^2(1-\delta)^2}
\sum_{j=1}^{M_{N,m,d}} \sum_{k=1}^n 
\left(
E\left[\Psi_j^2(X_k)g^2 (X_k)\right]
+\sigma^2 \cdot \mathbb E\left[\Psi_j^2(X_k)\right]
\right)\\
& = 
\frac{\left(h_0^{(\alpha)}\right)^{2d}}{n^2(1-\delta)^2}
\sum_{j=1}^{M_{N,m,d}} \sum_{k=1}^n 
\mathbb E\left[\Psi_j^2(X_k)g^2 (X_k)\right]
+
\frac{\left(h_0^{(\alpha)}\right)^{d} M_{N,m,d} \cdot \sigma^2}{n(1-\delta)^2}. 
    \end{aligned}
\end{equation}
Using Lemma~\ref{lm:Bounds}, we have :
$$ \mathbb E\left[\Psi_j^2(X_k)g^2 (X_k)\right]
\le \|\Psi_j||_{\infty}^2 \cdot \mathbb E\left[g^2 (X_k)\right]
\le  
\left(D(N,\alpha)\right)^{2m} \cdot
\frac{\|f-\Pi_{N,m}f \|_{\alpha}^2}{\left(h_0^{(\alpha)}\right)^{2d}}.
$$
This implies that
\begin{equation}
\begin{aligned}
\int_{\Omega_{+}} \| \Pi_{N,m}(f) - \widehat{f}_{N,n,m}^{(\alpha)}\|_{\alpha}^2d\rho_n & \le
   \mathbb E\left[ \| \Pi_{N,m}(f) - \widehat{f}_{N,n,m}^{(\alpha)}\|_{\alpha}^2
    \right]\\
&     \le 
     \frac{M_{N,m,d}}{n(1-\delta)^2}
     \bigg(
[D(N,\alpha)]^{2m} \cdot
\|f-\Pi_{N,m}f \|_{\alpha}^2
+\sigma^2 \left(h_0^{(\alpha)}\right)^{d} \bigg)\,.
\end{aligned}
\end{equation}
\qed

\section{Quality of the estimation in a weighted Sobolev space}

In this paragraph, we give a precise rate of convergence of our estimator in the case where the $d-$variate  regression functions belongs to a weighted Sobolev space. More precisely, we give an estimate for the term $\|f-\Pi_{N,m}f\|_{\alpha}^2$ in the $L^2-$risk error given by Theorem 2. For this purpose, we recall some definitions and results mainly borrowed from \cite{Potts}. 
Let 
$\mathbb T$ be a unit torus and let $f\in L^2(\mathbb T^d).$ The Fourier series expansion of $f$ is given by 
\begin{equation}\label{Eq5.1}
f(\pmb x)= \sum_{\pmb p\in \mathbb \Z^d} a_{\pmb p}(f) e^{2i\pi \pmb p\cdot \pmb x},\quad \pmb x\in \mathbb T^d,\quad (\mbox{or } \mathbb R^d).
\end{equation}
We have the following partition of $\Z^d$ : 
$\Z^d=\bigcup_{\pu \in \mathcal{D}} F_{\pu}$.
This gives the analysis of variance (ANOVA) decomposition of $f$ : 
$$ f(x)=f_{\emptyset}+\sum_{i=1}^d f_{i}(x_i)
+\sum_{i=1}^{d-1} \sum_{j=i+1}^d
f_{i,j}(x_i,x_j)+\ldots + f_{\mathcal{D}}(\px)
=\sum_{\pu \in \mathcal{D}} f_{\pu} (x_{\pu}) \,,$$ where the functions  $f_{\pu}$ are called ANOVA terms, see  \cite{Potts}.
Let $\pmb{U} \subset \mathcal{P}(\mathcal{D})$. The truncated ANOVA decomposition over $\pmb{U}$ is defined as:
$$ T_{\pmb{U}}(f):= \sum_{\pu \in \pmb{U}} f_{\pu} \,.$$
In particular, for $1 \le s \le d$, we define 
$$ T_s (f):=\sum_{|\pu| \le s } f_{\pu} \,.$$

\begin{definition}
Let $s>0$. Let $w^{(s)} : \Z^d \mapsto [1,\infty)$  be the weight function defined by
$$ w^{(s)}(\pmb{p}):= \prod_{j=1}^d \left(1 + |\pmb{p}_j|\right)^s\, \,  \  , \ \forall \ \pmb{p} \in \Z^d\,.$$
We associate to this weight function the Sobolev space 
\begin{equation}\label{weighted_Sobolev}
 H^{s}(I^d):=\left\{ f  \in L^2(I^d);\, 
\|f\|_{H^{s}(I^d)} :=\sum_{\pmb{p} \in \Z^d} (w^{(s)}(\pmb{p}))^2 \cdot |a_{\pmb{p}}(f)|^2 < +\infty \right\}\,,
\end{equation}
 and the weighted Wiener algebra
 $$\mathcal{A}^s(I^d):=\left\{ f  \in L^1(I^d);\, 
\|f\|_{\mathcal{A}^{s}(I^d)} = \sum_{\pmb{p} \in \Z^d} w^{(s)}(\pmb{p}) \cdot |a_{\pmb{p}}(f)| < +\infty \right\}\,. $$
\end{definition}

The following lemma provides us with an estimate for the decay rate of the expansion series expansion coefficients of the $f_{\pu}$ when $f\in  H^{s+\frac{d}{2}}(I^d)$ and with respect to the basis functions $\Phi_{\pk,\pmb{v}}^{(\alpha)}.$

\begin{lemma} \label{lm:boundC}
    Let $0 < \xi < \frac{1}{2}$, $N \in \N$ and  $\pk \in F_{\pu}$  such that  $||\pk||_1 \ge N+1$ with 
    $$\frac{\frac{N+1}{|\pu|}}{\ln\left(\frac{N+1}{|\pu|}\right)}
   \ge -\ln\left(\xi\right)
\left(1+s+\frac{d}{2} \right)\,.$$ For any 
$f \in H^{s+\frac{d}{2}}(I^d)$
\begin{equation}
    | C_{\pk,\pu}(f_{\pu})|=|< f_{\pu}, \Phi_{\pk,\pmb{v}}^{(\alpha)} >_{\alpha}| 
  \lesssim_{\alpha,d,s}
  || \pk||_{\infty}^{-s-\frac{d}{2}}
   \left(
    || f_{\pu}||_2 + ||f_{\pu}||_{H^{s+\frac{d}{2}}(I^d)}
    \right)\,.
\end{equation}
\end{lemma}
\proof~
We have:
\begin{equation}
f_{\pu}(\pmb x)= \sum_{\pmb p\in \mathbb \Z^d} a_{\pmb p}\left(f_{\pu}\right) e^{2i\pi \pmb p\cdot \pmb x},\quad \pmb x\in I^d.
\end{equation}
For $\pmb{v} \in \Z^d$ and $\pk \in \Z^{|\pmb{v}|}$, we have :
$$ C_{\pk,\pmb{v}}(f_{\pu})=
< f_{\pu}, \Phi_{\pk,\pmb{v}}^{(\alpha)} >_{\alpha}
=\sum_{\pmb p\in \mathbb \Z^d} a_{\pmb p}(f_{\pu})
< e^{2i\pi \pmb p\cdot \pmb x}, \Phi_{\pk,\pmb{v}}^{(\alpha)}(\px) >_{\alpha}
=
\sum_{\pmb p\in \mathbb \Z^d} a_{\pmb p}(f_{\pu}) d_{\pk,\pmb{v},\pmb{p}}\,,
$$
where $d_{\pk,\pmb{v},\pmb{p}}=< e^{2i\pi \pmb p\cdot \pmb x}, \Phi_{\pk,\pmb{v}}^{(\alpha)}(\px) >_{\alpha}$.
Note that  if $\pmb{v} \neq \pu$, then 
$$ C_{\pk,\pmb{v}}\left(f_{\pu}\right)=<f_{\pu}, \Phi_{\pk,\pmb{v}}^{(\alpha)}>_{\alpha}=0\,.$$
Consequently, we need only to estimate 
$ C_{\pk,\pu}\left(f_{\pu}\right)$.
In the special case $|\pu|=0$,  it is easy to show that $|d_p|:=|<e^{2i \pmb{p},\px}, \Phi_{\pmb{0}}^{(\alpha)}(\px)>_{\alpha}| \le 1$. Next, for  $|\pu|>0$,  we have
\begin{equation}
\label{eq:dbound}
    \begin{aligned}
    \big| d_{\pk,\pu,\pmb{p}} \big| & =
    \Big|<e^{2i \pmb{p},\px}, \left(\frac{1}{\sqrt{h_0^{(\alpha)}}}\right)^{d-|\pu|}
     \prod_{j \in \pu } \widetilde{P}_{k_j}^{(\alpha)}(x_j)
     >_{\alpha} \Big|  \\
      & = \Big|\int_{[-1,1]^d} \prod_{l=1}^d e^{2i \pi p_l x_l} 
      \left(\frac{1}{\sqrt{h_0^{(\alpha)}}}\right)^{d-|\pu|}
       \prod_{j \in \pu } \widetilde{P}_{k_j}^{(\alpha)}(x_j)
       \pmb{w}(\px)d\px  \Big|  \\
       &= \left[  \prod_{l \in \mathcal{D} \backslash \pu}
       \Bigg|  \int_{-1}^1 
      \frac{e^{2i \pi p_l x_l}}{\sqrt{h_0^{(\alpha)}}} w_{\alpha}(x_l)  dx_l \Bigg|
       \right] \times \left[\prod_{j \in  \pu} \Bigg| 
       \int_{-1}^1 
      e^{2i \pi p_j x_j}  \widetilde{P}_{k_j}^{(\alpha)}(x_j)  w_{\alpha}(x_j)  dx_j \Bigg|
       \right]\\
       & \le  \prod_{j \in  \pu} \Bigg| 
       \int_{-1}^1 
      e^{2i \pi p_j y}  \widetilde{P}_{k_j}^{(\alpha)}(y)  w_{\alpha}(y)  dy \Bigg|=\prod_{j \in  \pu} \Bigg| d_{p_j,k_j} \Bigg| \,,
    \end{aligned}
\end{equation}
where, for $l \in \Z$ and $r \in \N$, 
$d_{l,r}:=\int_{-1}^1       e^{2i \pi l y}  \widetilde{P}_{r}^{(\alpha)}(y)  w_{\alpha}(y)  dy$. Using (\ref{eq:Bessel1}), we get
$$ d_{l,r}=i^r \sqrt{\pi} \sqrt{2r+2\alpha+1}\sqrt{\frac{\Gamma(r+2\alpha+1)}{\Gamma(r+1)}}
\frac{J_{r+\alpha+\frac{1}{2}}(2\pi l ) }{(2\pi l )^{\alpha+\frac{1}{2}}}\,,
$$
where $J_a$ is the Bessel function of the first kind and order $a>-1$.
Note that since $d_{l,r}=(-1)^r d_{-l,r}$ then $|d_{l,r}|=|d_{-l,r}|$
which allow us to consider only the case where $l\ge 0$. Now, we apply the inequality (\ref{eq:Bessel2}) to the Bessel function $J_a$
of the previous equation. We obtain :
$$|d_{l,r}| \le  \sqrt{\pi} \sqrt{2r+2\alpha+1}\sqrt{\frac{\Gamma(r+2\alpha+1)}{\Gamma(r+1)}}
\frac{1}{\Gamma(r+\alpha+\frac{3}{2})}
\frac{(\pi |l| )^r}{2^{\alpha+\frac{1}{2}} }\,.
$$
Using (\ref{eq:Gamma}), we get
\begin{equation}
    \begin{aligned}
     |d_{l,r}|  & \le \sqrt{r} (e\pi |l|)^r   \frac{\pi^{\frac{3}{4}}e^{\frac{1}{4}}}{2^{\alpha+\frac{1}{2}}}
\sqrt{\alpha+\frac{3}{2}}
\frac{(r+2\alpha+\frac{1}{2})^{\frac{r}{2}+\alpha+\frac{1}{4}}}
{(r+\frac{1}{2})^{\frac{r}{2}+\frac{1}{4}} \cdot (r+\alpha+1)^{r+\alpha+1}}\\
& \le K(\alpha) \sqrt{r} \left( \frac{e \pi |l|}{r+\alpha} \right)^r\,, 
    \end{aligned}
\end{equation}
with $K(\alpha):=\frac{\pi^{\frac{3}{4}}e^{\alpha+\frac{1}{4}}}{\sqrt{2}}  \sqrt{\alpha+\frac{3}{2}}(\alpha+1)^{\alpha}$.
We also have, from Cauchy-Schwartz inequality that
$|d_{p_j,k_j}| \le 1$ for all $k_j \in \N$ and  $p_j \in \Z$. Injecting this in (\ref{eq:dbound}), we get
$$|d_{\pk,\pu,\pmb{p}}| \le |d_{p_{j_0}, k_{j_0}}|\,  $$
where $j_0 \in [[1,d]] $ is such that $||\pk||_{\infty}=|k_{j_0 }| $.
Let $0 < \xi < 1 $ such that 
$||\pmb{p}||_{\infty} \le \xi \frac{||\pk||_{\infty}}{e \pi}$. This implies that 
$$ |d_{\pk,\pu,\pmb{p}}| \lesssim_{\alpha} \sqrt{||\pk||_{\infty}}
\xi^{||\pk||_{\infty}} \ , \ \mbox{if} \ 
||\pmb{p}||_{\infty} \le \xi \frac{||\pk||_{\infty}}{e \pi}\,,
$$
where the notation $\lesssim_{\gamma}$  means in general  that the inequality is true up to constant depending only on a variable $\gamma$. Let us now re-write the expression of $C_{\pk,\pu}(f_{\pu})$ in the following manner :
$$C_{\pk,\pu}(f_{\pu})=
\underbrace{\sum_{ e||\pi||_{\infty} \le \xi ||\pk||_{\infty} } a_{\pmb p}(f_{\pu}) d_{\pk,\pu,\pmb{p}}}_{S_1}+
\underbrace{\sum_{ e||\pi||_{\infty} > \xi ||\pk||_{\infty} } a_{\pmb p}(f_{\pu}) d_{\pk,\pu,\pmb{p}}}_{S_2}.
$$
To get an upper bound for $|S_1|,$ we let 
 $\mathcal{A}:=\{\pmb{p} \in \Z^d : 
||\pmb{p}||_{\infty} \le \xi \frac{||\pk||_{\infty}}{e \pi} \}$. This set contains at most 
$\left( 
2\left[ \frac{\xi}{e \pi}|| \pk||_{\infty} \right]+1 
\right)^d$ elements where $[x]$ denotes the integer part of $x$. Hence, by using  Bessel's  and Cauchy-Schwartz inequalities, we obtain
\begin{equation}
\label{eq:S1Bound}
    \begin{aligned}
     |S_1|  = \Bigg| \sum_{\pmb{p} \in \mathcal{A}}  a_{\pmb p}(f_{\pu}) d_{\pk,\pu,\pmb{p}} \Bigg| & \le 
     \left( \left( \sum_{ \pmb{p} \in \mathcal{A}}  |a_{\pmb p}(f_{\pu})|^2\right) \cdot  
     \left( \sum_{\pmb{p} \in \mathcal{A}}   |d_{\pk,\pu,\pmb{p}}|^2 \right) \right)^{\frac{1}{2}}\\
     & \lesssim_{\alpha} \left(  
      \sum_{\pmb{p} \in \mathcal{A}} 
     ||\pk||_{\infty} \xi^{2 ||\pk||_{\infty}} 
     \right)^{\frac{1}{2}} \cdot || f_{\pu}||_2\\
     & \lesssim_{\alpha}  \sqrt{||\pk||_{\infty}} 
     \xi^{ ||\pk||_{\infty}} 
     \left( 
2\left[ \frac{\xi}{e \pi}|| \pk||_{\infty} \right]+1 
\right)^\frac{d}{2}|| f_{\pu}||_2\\
& \lesssim_{\alpha} || \pk||_{\infty}^{\frac{d+1}{2}}  \xi^{ ||\pk||_{\infty}} || f_{\pu}||_2\,.
    \end{aligned}
\end{equation}
Next, we get an upper bound for $|S_2|$. For thus purpose, we suppose that $f\in H^{s+\frac{d}{2}}(I^d)$. It follows, using Lemma3.9 of \cite{Potts}, 
that $f_{\pu} \in H^{s+\frac{d}{2}}(I^d)$ $\forall \,  |\pu| \le d$. 
Applying Cauchy-Schwartz inequality to $|S_2|^2$ followed by
 Bessel's inequality, we obtain
\begin{multline}
    |S_2|^2 \le
 \left( \sum_{ \pmb{p} \in \mathcal{A}^c}  
 |a_{\pmb p}(f_{\pu})|^2\right). 
     \left( \sum_{\pmb{p} \in \mathcal{A}^c}   |d_{\pk,\pu,\pmb{p}}|^2 \right) \\
     \le \left( \sum_{ \pmb{p} \in \mathcal{A}^c}  
 |a_{\pmb p}(f_{\pu})|^2\right)  \|\Phi_{\pu,\pk}^{(\alpha)} \pmb{w}^{(\alpha)} \|_2^2
 \le 
 \left( \sum_{ \pmb{p} \in \mathcal{A}^c}  
 |a_{\pmb p}(f_{\pu})|^2\right).||\Phi_{\pu,\pk}^{(\alpha)} \|_{\alpha}^2 \cdot \|\pmb{w}^{(\alpha)}\|_{\infty}^2
 \le 
 \left( \sum_{ \pmb{p} \in \mathcal{A}^c}  
 |a_{\pmb p}(f_{\pu})|^2\right)
 \,.
\end{multline}
Note that for $\pmb{p} \in \mathcal{A}^c$, we have $  
\left[\prod_{j=1}^d \left(1+|p_j| \right)^{s+\frac{d}{2}}
\right]^2 \ge ||\pmb{p}||_{\infty}^{2s+d} 
\ge
\left[ \xi \frac{||\pk||_{\infty}}{(e \pi)}\right]^{2s+d}
$. 
Consequently, for $f_u \in H^{s+\frac{d}{2}}(I^d) $, we have :
$$ |S_2|^2 \le 
\sum_{ \pmb{p} \in \mathcal{A}^c} |a_{\pmb p}(f_{\pu})|^2
\le
\left(\frac{e \pi }{\xi ||\pk||_{\infty}}\right)^{2s+d}
\sum_{ \pmb{p} \in \mathcal{A}^c} 
\left(1+\sum_{j=1}^d |p_j|^2\right)^{s+\frac{d}{2}} |a_{\pmb{p}}(f_{\pu})|^2
\le
\left(\frac{e \pi }{\xi ||\pk||_{\infty}}\right)^{2s+d}
||f_{\pu}||_{H^{s+\frac{d}{2}}(I^d)}^2\,.
$$
Hence, 
\begin{equation}
\label{eq:S2Bound}
    |S_2| \lesssim_{\alpha,d,s}  \| \pk\|_{\infty}^{-s-\frac{d}{2}}
    \|f_{\pu}\|_{H^{s+\frac{d}{2}}(I^d)}\,.
\end{equation}
Combining (\ref{eq:S1Bound}) and (\ref{eq:S2Bound}), we get
\begin{equation}
    \begin{aligned}
     |C_{\pk,\pu}(f_{\pu})| 
&  \lesssim_{\alpha,d,s}
\| \pk\|_{\infty}^{\frac{d+1}{2}}  \xi^{ \|\pk||_{\infty}} \| f_{\pu}\|_2 + \| \pk\|_{\infty}^{-s-\frac{d}{2}}
    \|f_{\pu}\|_{H^{s+\frac{d}{2}}(I^d)}\\
&    \le
    \left( \| \pk\|_{\infty}^{\frac{d+1}{2}}  \xi^{ \|\pk\|_{\infty}}
    + \| \pk\|_{\infty}^{-s-\frac{d}{2}}
    \right)
    \left(
    \| f_{\pu}\|_2 + \|f_{\pu}\|_{H^{s+\frac{d}{2}}(I^d)}
    \right)
    \end{aligned}
\end{equation}
Finally, if $\|\pk\|_1 \ge (N+1)$ and $N \in \N$ is such that $\frac{N+1}{|\pu|.(\ln(N+1)-\ln(|\pu|))} \ge \frac{1}{\ln\left(\frac{1}{\xi}\right)}
\left(1+s+\frac{d}{2} \right),$ then\\
$\| \pk\|_{\infty}^{\frac{d+1}{2}}  \xi^{ \|\pk\|_{\infty}}
    \le \| \pk\|_{\infty}^{-s-\frac{d}{2}}$. Consequently
$$  |C_{\pk,\pu}(f_{\pu})| 
  \lesssim_{\alpha,d,s}
  \| \pk \|_{\infty}^{-s-\frac{d}{2}}
   \left(
    \| f_{\pu}\|_2 + \|f_{\pu}\|_{H^{s+\frac{d}{2}}(I^d)}
    \right)\,.
  $$
\qed

\begin{theorem}
\label{th:quality}~
Let $f \in H^{s+\frac{d}{2}}(I^d)$ and let  $ 1 \le m \le \min \{N,d\}$ with $N >e.m-1$ and
  such that 
  $$
  \frac{\frac{N+1}{m}}{\ln\left(\frac{N+1}{m}\right)}
   \ge -\ln\left(\xi\right)
\left(1+s+\frac{d}{2} \right)\,,$$
where $0 < \xi < \frac{1}{2}$. Suppose that $s+\frac{d}{2} > m+1,$ then
\begin{equation}\label{Estimate1}
    \|f-\Pi_{N,m}(f)\|_{\alpha} \lesssim_{\alpha,d,s} \frac{1}{2^{(s+\frac{d}{2})(m+1)}}
    ||f||_{H^{s+\frac{d}{2}}(I^d)} 
     + \frac{m \|f\|}{s+\frac{d}{2}-m}. \left(2+\frac{m}{N+1}\right)^{d} \left(\frac{m}{N+1}\right)^{s+\frac{d}{2}-m}    \,.
\end{equation}
and
\begin{equation} 
\|f-\Pi_{N,m}(f)\|_{\infty} \lesssim_{\alpha,d,s} 
     \frac{1}{2^{(s+\frac{d}{2})(m+1)}}\left\|f \right\|_{\mathcal{A}^s(I^d)}+
       \frac{\|f \| D^m(N,\alpha)m}{s+\frac{d}{2}-m}
        \left(2+\frac{m}{N+1}\right)^{d}
        \left( \frac{m}{N+1}\right)^{s+\frac{d}{2}-m}  \,.
\end{equation}
Here, $\|f\|=||f||_2+||f||_{H^{s+\frac{d}{2}}}$ and 
$D(N,\alpha)$ is as given by Lemma~\ref{lm:Bounds}.
\end{theorem}
\proof~ The orthogonal projection of $f$ on $\mathcal{P}_{N,m,d}$, $\Pi_{N,m} f$, verifies
$$\left\|f- \Pi_{N,m} f\right\|_\alpha =\Big\| 
\big[ f-T_m(f)\big] 
+\big[ T_m(f) - \Pi_{N,m}(T_m(f)) \big]
+\big[  \Pi_{N,m}(T_m(f)) -\Pi_{N,m}(f)\big]
\Big\|_\alpha $$
Note that for $N \ge m$, $\Pi_{N,m}(T_m(f)) =\Pi_{N,m}(f)$ and consequently, we have
$$\left\|f- \Pi_{N,m} f\right\|_\alpha \le
\left\|  f-T_m(f)\right\|_\alpha 
+\left\| T_m(f) - \Pi_{N.m}(T_m(f)) \right\|_\alpha. $$
From \cite{Potts}, we have for $f \in H^{s+\frac{d}{2}}(I^d)$
$$ \| f- T_mf \|_{L^2(I^d)} \le \frac{1}{2^{(s+\frac{d}{2})(m+1)}} \|f\|_{H^{s+\frac{d}{2}}(I^d)}\,.$$
To get an upper bound for $ \left\| T_m(f) - \Pi_{N,m}(T_m(f)) \right\|_\alpha$, we proceed as follows.  
The function $\Pi_{N,m} f$ writes as 
$ \Pi_{N,m} f=\sum_{|\pmb{v}| \le m ; ||\pk ||_1 \le N }
C_{\pk,\pmb{v}}(f) \Phi_{\pk,\pmb{v}}(x)\,.$
hence, 
\begin{equation}
\begin{aligned}
 \left\| T_m(f) - \Pi_{N,m}(T_m(f)) \right\|_{\alpha} & =
\Big\| \sum_{|\pu| \le m} f_{\pu} - \Pi_{N,m}\Big(\sum_{|\pu| \le m} f_{\pu}\Big) \Big\|_{\alpha} \\
& =
\left\| \sum_{|\pu| \le m} \left[f_{\pu} - \Pi_{N,m}( f_{\pu})\right] \right\|_{\alpha} 
 = \left\| \sum_{|\pu|\le m ; |\pmb{v}| \le m ; \|\pk\|_1 \ge N+1}
C_{\pk,\pmb{v}}\left(f_{\pu}\right) \Phi_{\pk,\pmb{v}}^{(\alpha)}
\right\|_{\alpha} \\
&= \left\|  \sum_{|\pmb{u}| \le m ; \|\pk\|_1 \ge N+1}
C_{\pk,\pu}\left(f_{\pu}\right) \Phi_{\pk,\pu}^{(\alpha)}
\right\|_{\alpha} \,.
\end{aligned}    
\end{equation}
Consequently,  we get
\begin{equation}
    \begin{aligned}
     \left\| T_m(f) - \Pi_{N,m}(T_m(f)) \right\|_{\alpha} & \le 
       \sum_{|\pmb{u}| \le m ; \|\pk\|_1 \ge N+1}
     |C_{\pk,\pu}\left(f_{\pu}\right)| \cdot \left\|\Phi_{\pk,\pu}^{(\alpha)}\right\|_{\alpha}=\sum_{|\pmb{u}| \le m ; \|\pk_{\pu}\|_1 \ge N+1}
     |C_{\pk,\pu}\left(f_{\pu}\right)|\\
     & = 
      \sum_{|\pmb{u}|=1}^m \sum_{\|\pk\|_1 \ge N+1}
     |C_{\pk,\pu}\left(f_{\pu}\right)|\\
     & \le 
      \sum_{|\pmb{u}|=1}^m \sum_{\|\pk\|_{\infty} \ge \frac{N+1}{|\pu|}}
     |C_{\pk,\pu}\left(f_{\pu}\right)|\,.
    \end{aligned}
\end{equation}
Let $1 \le |\pu| \le m$ and $\pk \in F_{\pu}.$ Then 
for $1 \le j \le m$ and $l \in \N$, the number of elements $(\pu,\pk)$ such that   
$|\pu|=j$ and
$||\pk ||_{\infty}=l$
 is ${d \choose j}
j (l+1)^{j-1}\,.$ 
Let $\|f_{\pu}\|:=||f_{\pu}||_2+||f_{\pu}||_{H^{s+\frac{d}{2}}}$. Using the result of Lemma3.9 in \cite{Potts} and adapting its proof for the $||.||_{H^{s}}-$norm, we conclude  that for all $\pu \subset \mathcal{D}$ with $|\pu| \ge 1$
 $\|f_{\pu}\| \le  \|f\|.$ 
Moreover, by using  the result of Lemma~\ref{lm:boundC}, we get 
\begin{equation}
    \begin{aligned}
     \left\| T_m(f) - \Pi_{N,m}(T_m(f)) \right\|_{\alpha} &  \lesssim_{\alpha,d,s}  
      \|f\| \cdot \sum_{j=1}^m 
      {d \choose j} \, j \, 
       \sum_{l \ge \frac{N+1}{j}} (l+1)^{j-1} \frac{1}{l^{s+\frac{d}{2}}}\\
       & \le 
     \|f\| \cdot \sum_{j=1}^m 
      {d \choose j} \,  j\,   \left(1+\frac{m}{N+1}\right)^{j-1}
       \sum_{l \ge \frac{N+1}{j}}  l^{j-1-s-\frac{d}{2}}\\
        & \le 
         \frac{m \|f\|}{s+\frac{d}{2}-m}. \left(2+\frac{m}{N+1}\right)^{d} \left(\frac{m}{N+1}\right)^{s+\frac{d}{2}-m}\,.
    \end{aligned}
\end{equation}
In a similar manner,   we have
$$\left\|f- \Pi_{N,m} f\right\|_{\infty} \le
\left\|  f-T_m(f)\right\|_{\infty} 
+\left\| T_m(f) - \Pi_{N.m}(T_m(f)) \right\|_{\infty} \,.$$
In  \cite{Potts}, it has been shown that 
$$\left\|  f-T_m(f)\right\|_{\infty}  \le  \frac{1}{2^{(s+\frac{d}{2})(m+1)}}\left\|f \right\|_{\mathcal{A}^s(I^d)}. $$
On  the other hand and by using Lemma~\ref{lm:Bounds} and Lemma~\ref{lm:boundC}, one gets
\begin{eqnarray*}
     \left\| T_m(f) - \Pi_{N.m}(T_m(f)) \right\|_{\infty} & \le&
       \sum_{|\pmb{u}|=1}^m \sum_{\|\pk\|_{\infty} \ge \frac{N+1}{|\pu|}}
     |C_{\pk,\pu}\left(f_{\pu}\right)|. \left\|\Phi_{\pk,\pu}^{(\alpha)}\right\|_{\infty}\\
     & \le&
       \sum_{|\pmb{u}|=1}^m \sum_{\|\pk\|_{\infty} \ge \frac{N+1}{|\pu|}} |C_{\pk,\pu}\left(f_{\pu}\right)|.
        \cdot \frac{D^m(N,\alpha)}{
        \left(h_0^{(\alpha)}\right)^{\frac{d}{2}}}\\
       & \lesssim_{\alpha,d,s} & \frac{D^m(N,\alpha)}{
        \left(h_0^{(\alpha)}\right)^{\frac{d}{2}}} \cdot
       \|f\| 
       \sum_{|\pmb{u}|=1}^m \sum_{\|\pk\|_{\infty} \ge \frac{N+1}{|\pu|}}
       ||\pk||^{-s-\frac{d}{2}}_{\infty}  \\
       & = &
       \frac{D^m(N,\alpha)}{
        \left(h_0^{(\alpha)}\right)^{\frac{d}{2}}}
        \|f\| 
       \sum_{j=1}^m 
       {d \choose j} \,  j\,  
       \sum_{l \ge \frac{N+1}{j}} 
       l^{j-1-s-\frac{d}{2}} \left(1+\frac{1}{l}\right)^{j-1} \\
       & \lesssim_{\alpha,d,s} &
        \|f\|  D^m(N,\alpha)
       \sum_{j=1}^m      {d \choose j} \, j\,  \left(1+\frac{j}{N+1}\right)^{j-1}
       \sum_{l \ge \frac{N+1}{j}} 
       l^{m-1-s-\frac{d}{2}} \\
       & \lesssim_{\alpha,d,s} &
           \frac{m \|f\|  D^m(N,\alpha)}{s+\frac{d}{2}-m}
        \left(2+\frac{m}{N+1}\right)^{d}
        \left( \frac{m}{N+1}\right)^{s+\frac{d}{2}-m}\,. 
\end{eqnarray*}
\qed

\begin{remark} From \eqref{Estimate1}, one concludes that under the conditions that ${\displaystyle \gamma_{m,N}=\frac{m}{N+1}<1}$ and $s+\frac{d}{2}-m>0$ with $2^{(s+\frac{d}{2})(m+1)} > 3^d \gamma_{m,N}^{m-\frac{d}{2}-s},$ we have 
$$  \|f-\Pi_{N,m}(f)\|_{\alpha} = O\Big( 3^d \gamma_{m,N}^{s+\frac{d}{2}-m}\Big).$$
The previous estimate together with the estimate given by Theorem 2, provide us with a precise estimate for the $L^2-$risk error of our estimator 
$\widehat F^{(\alpha)}_{N,n,m}$ when the regression function belongs to the weighted Sobolev space \eqref{weighted_Sobolev}.
\end{remark}

\section{Numerical Examples}\label{num_results}

In this paragraph, we give three numerical examples that illustrate the different results of this work.\\

\noindent
{\bf Example 1:} In this first example, we illustrate the results of Proposition~\ref{lm:DimOrder} and Theorem~\ref{th:CondNumber}. For this purpose, we have considered the following parameter values:
$d=4$ and $d=6,$ for the dimension  with the values of $2\leq N\leq 5,$ for the parameter relative to the total degree $d-$variate Jacobi polynomials space. These polynomials are associated with the two special values of $\alpha=-0.5$ and $\alpha=0.5.$
Moreover, we have considered the value of $m=2,$ that we restrict ourselves to  the ANOVA decomposition with $m=2$  interactions between the covariables. Also, for the case of $d=4,$ (respectively $d=6$), we have considered a random sampling set with size $n=900$
(respectively $n=1600$) and following a multivariate   $\mbox{Beta}(\alpha+1,\alpha+1)$ distribution. Then, we have computed the  $\kappa_2\left(G_{N,n,m}^{(\alpha)}\right),$ the $2-$condition number of the random projection matrix  $G_{N,n,m}^{(\alpha)},$ given by \eqref{eq:randomamatrix}. Also, for each values of the couple $(d,N),$ we have provided the dimension $M_{N,2,d}$ of the considered $d-$variate polynomial space. The obtained numerical results are given by the following Table~\ref{tab1}. 

\begin{center}
\begin{table}[h]
\vskip 0.2cm\hspace*{1.0cm}
\begin{tabular}{cccccccc} \hline
 $d=4$   &$M_{N,2,d}$&$\kappa_2
 \left(G_{N,n,2}^{(-0.5)}\right)$& $\kappa_2(G_{N,n,2}^{(0.5)})$& $d=6$   &$M_{N,2,d}$&$\kappa_2\left(G_{N,n,2}^{(-0.5)}\right)$& $\kappa_2\left(G_{N,n,2}^{(0.5)}\right)$ \\ 
 $N$& & & &$N$& & &  \\
 & & &  & & & & \\  \hline
  $2$& $15$ & $8.80$  &  $6.70$&$2$& $28$ & $16.90$  &  $11.75$ \\
  $3$& $31$ & $20.32$  & $32.50$&$3$& $64$ & $24.35$  &  $47.80$ \\
  $4$& $53$ & $41.16$  & $48.25$&$4$& $115$ & $54.20$  &  $87.75$ \\
  $5$& $81$ & $65.19$  & $135.45$&$5$& $181$ & $223.22$  &  $312.77$ \\ 
 \end{tabular}
\caption{The Jacobi polynomial space
dimension
and the $2-$condition number of the random projection matrix for $m=2,$ $(d,n)=(4,900),\, (6,1600),$
$\alpha=\frac{1}{2}, -\frac{1}{2}$ and $2\leq N\leq 5.$}
\label{tab1}
\end{table}
\end{center}
 
 Also, in  Figure~1, we give the plots of the spectrum of $G_{N,n,2}^{(-0.5)},$ for $(d,n)=(4,900), (6,1600)$ and the different values of $2\leq N\leq 5.$ Note that these plots are fairly coherent with the predicted theoretical behaviour of
the spectrum of the random matrix $G_{N,n,2}^{(\alpha)},$
given by Theorem~\ref{th:CondNumber} and  in particular by the lower and upper bounds \eqref{eq:lambdamin} and \eqref{eq:lambdamax}.

\begin{center}
\begin{figure}
\includegraphics[scale=0.7]{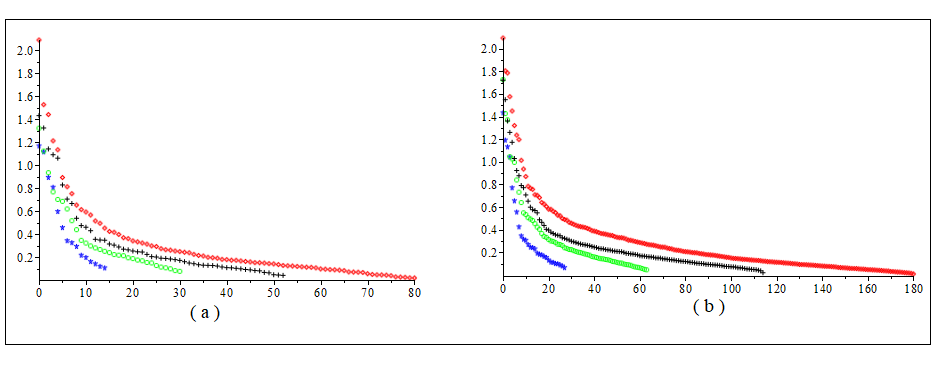}
\caption{Plots of the spectrum of $G_{N,n,2}^{(-0.5)},$ for (a) $(d,n)=(4,900)$ (b) $(d,n)=(6,1600)$ and  different values $N=2,3,4,5,$ (from left to right).}
\label{fig2}
\end{figure}
\end{center}

\noindent
{\bf Example 2:} In this second example, we illustrate the performance of $\widehat{f}_{N,n,m}^{(\alpha)},$ our  proposed stable NP regression  estimator, that is based on  least squares by means of  multivariate Jacobi polynomials. For this purpose, we have considered the NP regression problem \eqref{eq:Pb} for the special case of the dimension $d=4$ with a synthetic test true regression function $f$, given by 

\begin{eqnarray}\label{testfct1}
f(x,y,z,t)&=& x+(2 y-1)^2+\frac{\sin(2 \pi z)}{2-\sin(2 \pi z)}+0.1\, \sin(2 \pi t) +0.2 \, \cos(2 \pi t) \nonumber \\
&& \hspace*{5cm}+0.3\, (\sin(2 \pi t))^2+0.4\, (\cos(2 \pi t))^3+0.5\,  (\sin(2 \pi t))^3. 
\end{eqnarray}
Note that  this test regression function corresponds to an additive multidimensional regression model. Hence,  $m=1$ is the appropriate value of this parameter. Then, we have constructed our estimator 
$\widehat{f}_{N,n,1}^{(\alpha)}$ with $\alpha=-\frac{1}{2},$ $n=900$ i.i.d. random sampling points following a $4-$D $\mbox{Beta}(\alpha+1,\alpha+1)$ distribution. Also, we have considered the different values of $N=4,6,8, 10$ together with a noise free model as well as  noised models associated to two different values of $\sigma=0.1, 0.5.$ We have computed the empirical mean squared error over a test random set of size $n$ with i.i.d. random points $X_i$ following also a 
$4-$D $\mbox{Beta}(\alpha+1,\alpha+1)$ distribution.
This empirical mean squared error is given by 
$$MSE = \frac{1}{n} \sum_{i=1}^n 
\Big(\widehat{f}^{(\alpha)}_{N,n,1}(X_i)-f(X_i)\Big)^2.$$
The obtained numerical results are given by the following Table~\ref{Tab2}.

\begin{center}
\begin{table}[h]
\vskip 0.2cm\hspace*{5.0cm}
\begin{tabular}{cccc} \hline
   &$\sigma=0$&$\sigma=0.1$& $\sigma=0.5$\\  $N$&$MSE $&$MSE$& $MSE$ \\ 
   & & & \\ \hline
   & & & \\
 $4$ & $2.69 e-2$ & $3.99 e-2$ & $ 4.42 e-2$ \\
 $6$ & $1.60 e-2$ & $2.11 e-2$ & $ 2.53 e-2$ \\
  $8$ & $5.89 e-3$ & $1.05 e-2$ & $ 1.79 e-2$ \\
  $10$ & $9.29 e-4$ & $2.63 e-3$ & $ 6.20 e-2$ \\
  & & &   \\ \hline
  \end{tabular}
\caption{Numerical simulations for test function~(\ref{testfct1}).}
\label{Tab2}
\end{table}
\end{center}
Note that the numerical results given by Table~\ref{Tab2} are coherent with  the theoretical $L^2-$risk and the approximation error of the proposed estimator $\widehat{f}_{N,n,1}^{(\alpha)},$ given by Theorems~\ref{th:l2risk} and~\ref{th:quality}. Also, note that  the loss of accuracy we have observed for the special values of  $N=10$ and $\sigma=0.5$ is due to the fact that for the given couple $(N,n)=(10,900),$ the $2-$condition number of the random matrix $G_{N,n,m}^{(\alpha)}$ is relatively large to handle
noised data with relatively large $\sigma=0.5$ Nonetheless,  for $n=1600$ and  the same values of the parameters, we have obtained an $MSE \approx 1.16 e-2.$ \\

\noindent
{\bf Example 3:} In this last example, we consider the Kriging model test function borrowed from  \cite{Chen} and given, for $x_1,x_2,x_3,x_4 \in [0,1] $,  by 
\begin{equation}\label{testfct2}
f(x_1,x_2,x_3,x_4)= 1+ \exp\left[-2 \big((x_1-1)^2+x_2^2\big)-0.5 \big(x_3^2+x_4^2\big)\right]+ \exp\left[-2 \big(x_1^2+(x_2-1)^2\big)-0.5 \big(x_3^2+x_4^2\big)\right]\,. 
\end{equation}
Then, we have considered the values of the parameters $m=2,$
$\alpha=-\frac{1}{2},$ $n=1600$ and constructed our estimator
$\widehat{f}_{N,n,2}^{(\alpha)},$ with $N=4,5,6.$ As in the previous example, we have computed the different $MSE:$ the empirical mean squared errors for the different values of the Gaussian noise standard deviation $\sigma=0,\, 0.1,\,  0.5.$ These $MSE$ are computed by the use of a new set of $n_1=400$ i.i.d. random sampling points following a multivariate $\mbox{Beta}(\alpha+1,\alpha+1)$ distribution. The obtained numerical results are given by Table~\ref{Tab3}.

\begin{center}
\begin{table}[h]
\vskip 0.2cm\hspace*{5.0cm}
\begin{tabular}{cccc} \hline
   &$\sigma=0$&$\sigma=0.1$& $\sigma=0.5$\\  $N$&$MSE $&$MSE$& $MSE$ \\ 
   & & & \\\hline
   & & & \\
 $4$ & $5.35 e-3$ & $5.94 e-3$ & $ 1.38 e-2$ \\
 $5$ & $1.57 e-3$ & $1.52 e-2$ & $ 2.57 e-2$ \\
  $6$ & $1.06 e-3$ & $1.07 e-2$ & $ 2.69 e-2$ \\
   & & &   \\ \hline
  \end{tabular}
\caption{Numerical simulations for test function~(\ref{testfct2}) .}
\label{Tab3}
\end{table}
\end{center}
Note that these numerical results are also coherent with  the theoretical $L^2-$risk and the approximation error  given by Theorems~\ref{th:l2risk} and~\ref{th:quality}. For  the noise free model, that is $\sigma=0,$ the larger $N,$ the smallest is the empirical mean squared error. For the noised versions of the model, that is for $\sigma =0.1$ or $\sigma=0.5,$ the situation is slightly reversed. This is due to the contribution of the variance term. According to Theorem~\ref{th:l2risk}, this last quantity  is affected by larger values of the parameter $N.$




\end{document}